\documentclass{article}

\usepackage{amsfonts, amsmath, amssymb, amsthm, array, booktabs, cancel, cite, color, enumerate, graphbox, graphics, graphicx, hyperref, import, mathrsfs, mathtools, pdfpages, url, tabularx, tikz-cd, tikz,tikz-3dplot,,url, varwidth, wrapfig}
\usepackage[shortlabels]{enumitem}
\usepackage[all]{xy}
\usepackage[nottoc]{tocbibind}
\usepackage[margin=10pt,font=small,labelsep=period]{caption}

\usepackage[normalem]{ulem} 

\newlength{\notewidth}
\newtheorem{lemma}{Lemma}[section]
\newtheorem{prop}[lemma]{Proposition}
\newtheorem{example}[lemma]{Example}
\newtheorem{thm}{Theorem}
\newtheorem{cor}[lemma]{Corollary}

\theoremstyle{definition}
\newtheorem{rmrk}[lemma]{Remark}

\newcommand{\R}{{\mathbb R}}

\newcommand{\g}{\gamma}

\newcommand{\be}{\begin{equation}}
\newcommand{\ee}{\end{equation}}
\newcommand{\bp}{\begin{proof}}
\newcommand{\ep}{\end{proof}}

\newcommand{\cca}{\,{\bf c}\alpha_1 \,}
\newcommand{\ccb}{\,{\bf c}\alpha_2 \,}
\newcommand{\ssa}{\,{\bf s}\alpha_1 \,}
\newcommand{\ssb}{\,{\bf s}\alpha_2 \,}

\title{Geometry of Integrable Linkages}

\author{
Ron Perline
\footnote{
Department of Mathematics,
Drexel University
Philadelphia, PA 19104;
perlinrk@drexel.edu
}
\and
Serge Tabachnikov\footnote{
Department of Mathematics,
Penn State University, 
University Park, PA 16802;
tabachni@math.psu.edu}
}

\begin{document}

\maketitle

\begin{abstract}
In analogy with the well-known 2-linkage tractor-trailer problem, we define a 2-linkage problem  in the plane
with  novel non-holonomic ``no-slip''
conditions. Using constructs from sub-Riemannian geometry, we look for geodesics corresponding to 
linkage motion with these constraints (``tricycle kinematics'').  The paths of the three vertices turn out to
be critical points for functionals which appear in the hierarchy of conserved quantities for the planar
filament equation, a well known completely integrable evolution equation for planar curves.  We
show that the geodesic equations are completely integrable, and  present a second connection
to the planar filament equation.
\end{abstract}

\tableofcontents

\section{Introduction} \label{sect:intro}


There has been a good deal of interest recently in what is sometimes informally called ``bicycle mathematics".  This refers to questions about the geometry of  the kinematics of directed line segments in some space (usually $\R^2$ or $\R^3$) subject to a  constraint, the ``no-slip" condition, which we now describe. 

Given a  moving oriented line segment of a fixed length, we demand that  the velocity of the rear point be  parallel to the segment.  One imagines the front point to be the front wheel of the bicycle, the second point corresponds to the rear wheel. There are a great number of surprises in the subject, which has a long history and a recent renaissance; for background and history, we refer to \cite{BLPT}, \cite{FLT}, and \cite{Ta} for more information.

Our focus is on a recent development in the subject involving a fundamental minimization problem:  given two  oriented line segments $S_1, S_2$ of equal length, we wish to move $S_1$ to the position of $S_2$, where the motion is constrained by the no-slip condition indicated above.  How can we move $S_1$ so that the length of the path of the front vertex is minimized?  What are the ``bicycling geodesics"?

As far as we know, this question was first posed by Maxim Arnold, and addressed in \cite{2D}.  In that paper, the authors use the language and constructs of sub-Riemannian geometry to define a Hamiltonian on the cotangent bundle to the space of segments; the trajectories of the Hamiltonian flow project to arc length-parameterized curves which, at least locally, give the trajectories of shortest distance.  

 Remarkably, it is shown in \cite{2D} that these planar curves are {\it elasticae} - critical curves for the functional $\int_\gamma (\kappa^2  + \lambda) \, ds$, where  $\kappa$ is the curvature and the integration is with respect to arc length.  We will refer to $\int_\gamma \kappa^2  \, ds$ as the {\it elastic (or bending) energy}.  
 
 An analogue of this result exists in $\R^3$ also:  if a segment moves through space subject to the no-slip condition, then it is a Kirchhoff elastic rod, a critical point of $\int_\gamma  (\kappa^2  + \mu \tau+ \lambda) \, ds$, where $\tau$ is the torsion of the space curve. The bicycling geodesics in Euclidean spaces of higher dimensions are confined to affine 3-dimensional spaces, and the problem reduces to $\R^3$, see \cite{BJT}.

The authors of \cite{2D} give examples of the mysterious appearance of elastica in many geometric contexts, in addition to the one discovered in their paper.  There are more \cite{P}, and a phenomenon we have observed is that the presence of elastica in a geometrical context often indicates that there is an integrable system floating around in the background. 

 We now discuss one such context, which will be important to us.
Consider the evolution equation on planar curves, the {\it planar filament flow}, see \cite{LP},
$$\gamma_t = \frac{\kappa^2}{2}T + \kappa_s  N,$$
where $(T,N)$ is the Frenet frame, and $\kappa$ is the curvature.

This evolution on curves is locally arc-length preserving, so the curve bends without stretching or shrinking.  If the curve is closed, then the  total length $I_0=\int_\gamma 1 \, ds$ is conserved.  In addition, one can check that the elastic energy is also conserved. 

 In fact, there is an infinite list of such conserved quantities, whose existence can be explained by the close correspondence between the planar filament equation and the modified KdV equation:
$$
\kappa_t = \kappa_{sss}  + \frac{3}{2}\kappa^2 \kappa_s.
$$
If a curve evolves via the planar filament equation, then its curvature evolves via mKdV.
The mKdV equation, like its better known ``cousin" the KdV  equation, is completely integrable, and the planar filament equation inherits this integrability. 

  Among other things, this means that there is an infinite list of conserved quantities, defined in terms  of the curvature and its derivatives.  The elastic energy is one of these conserved quantities.  We list
the first three, as well as their Euler-Lagrange operators: 
\begin{equation*}
\begin{aligned}
I_0&=\int_\gamma 1 \, ds ,  \quad E_0 = \kappa(s) ,\\
I_2&=\int_\gamma \kappa^2  \, ds, \quad E_2 = \kappa'' + \frac{\kappa^3}{2},\\
I_4&=\int_\gamma (\kappa')^2- \frac{1}{4}\kappa^4 \, ds, \quad  E_4  = 
\frac{5}{2} \kappa^{2} \mathit{\kappa''} +\frac{5}{2} \kappa  \,\mathit{\kappa'}^{2}+\frac{3}{8} \kappa^{5}
+\mathit{\kappa''''}.
\end{aligned}
\end{equation*}
We will refer to critical points  of linear combinations of these functionals as {\it soliton curves}; for example, a curve satisfying $a E_0 + E_2 = 0$ is a 1-soliton, and a curve satisfying $a E_0 + b E_2 + E_4=0$ is a 2-soliton.

With this background in mind, we look at the appearance of elastica in \cite{2D}  in the study of  a single line segment, and ask if there might be similar interesting phenomena related to a pair of linked segments. 

\begin{figure}[ht]
\centering
\includegraphics[width=.4\textwidth]{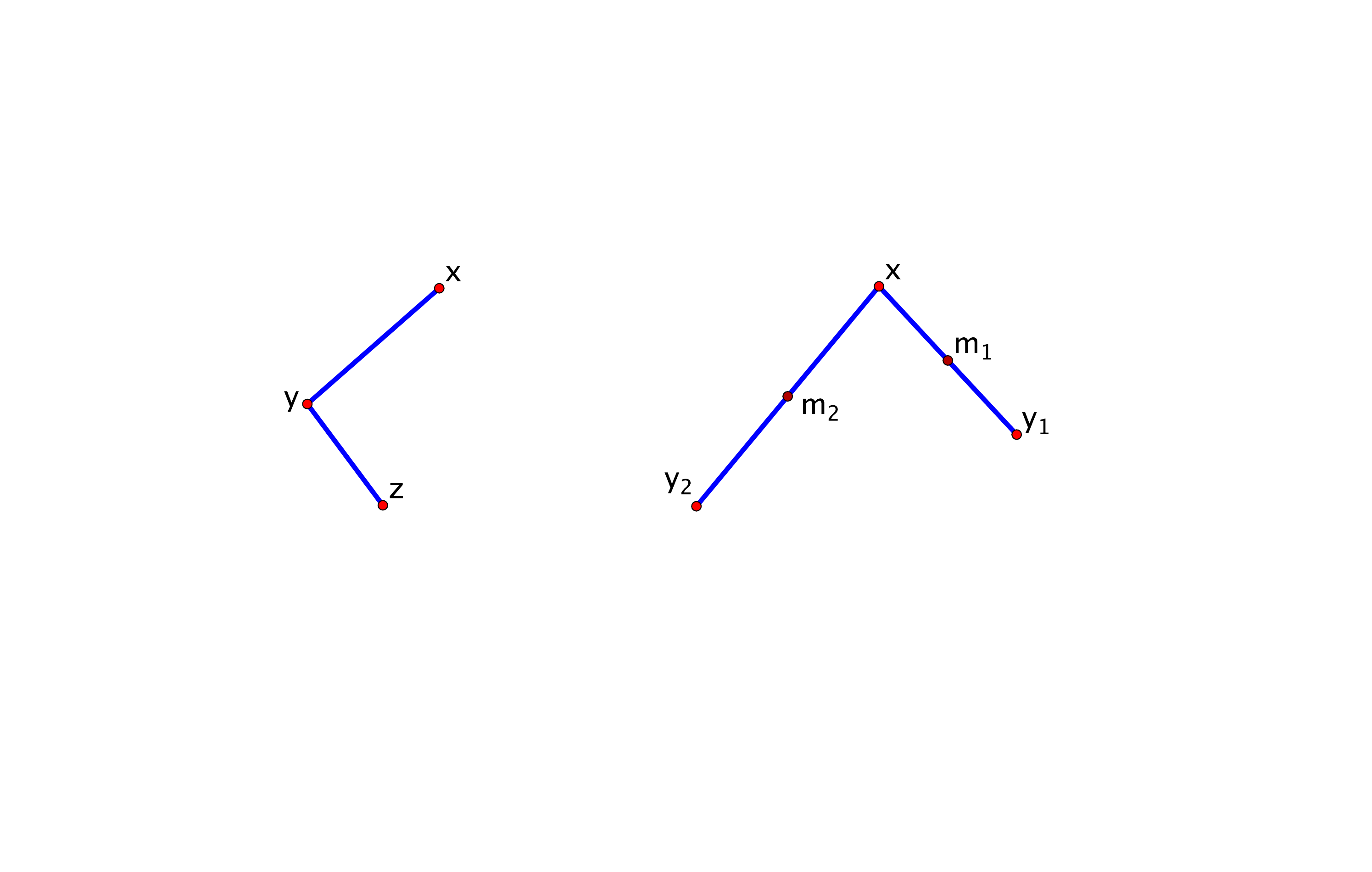}
\caption{Two-linkage: the no-slip condition is imposed at points $m_1$ and $m_2$.}	
\label{link}
\end{figure}

In analogy with the single segment, we seek to impose appropriate non-holonomic constraints for the motion, construct its associated Hamiltonian, and consider the geodesics on 2-linkage space in the plane, depicted in Figure \ref{link}.  

The segments $xm_1$ and $xm_2$ have fixed lengths (not necessarily equal). 
The no slip condition is imposed on points $m_1$ and $m_2$, whereas point $x$ can move without restriction. One may think of the linkage $x m_1 m_2$ as a kind of a {\it tricycle}, where $x$ is the front wheel and $m_{1}$ and $m_2$ are the two back wheels. 

The segments $xm_1$ and $xm_2$ are extended to $xy_1$ and $xy_2$ by doubling their lengths. As we will see, the no slip condition implies that points $y_1$ and $y_2$ will move with the same speed as point $x$ and therefore their trajectories have the same lengths as that of point $x$.

The variational problem  is to describe the {\it tricycling geodesics}, the motions of the linkage subject to the above no slip constraint that locally minimize the length of the trajectory of the front point $x$ (and hence those of points $y_!$ and $y_2$ as well). In particular, we want to describe the trajectories of the points $x,y_1$, and $y_2$, the projections of the respective sub-Riemannian geodesics to the plane.

We now turn to the contents of the paper.  In Section \ref{sect:1link}, following \cite{2D}, we review the equations of 1-linkage geodesic. 
In Section \ref{sect:distr} we show that the distribution in the configuration space of the 2-linkage, described by the no-slip constraint, is completely non-integrable, and in Section \ref{sect:sing}
we describe the singular curves of this distribution.  

In Section \ref{sect:eq}, we write down the equations of geodesics for 2-linkages and show how 1-linkage geodesics ``lift" to 2-linkage geodesics.
Next, in Section \ref{sect:ela}, we study the problem of 2-linkages when the segments have equal length.  Surprisingly, it turns out that the path of the ``front" point $x$ of the 2-linkage is a 1-soliton (elastica). Unlike the case of 1-linkages, we obtain both non-inflectional and inflectional elasticae. 

Section \ref{sect:comp} presents a different, computer-assisted, method of obtaining the results of Section \ref{sect:ela}, using Gr\"obner bases. We use this method later in the paper, when the calculations become too overwhelming.

Section \ref{sect:Ba} is an aside on using B\"acklund transformations (defined at that point) to take 1-soliton curves to 2-soliton curves.  Using these results, we conclude  that (in the equal lengths case) the paths associated with the other two end points are 2-soliton curves.

Section  \ref{sect:examples} presents a variety of examples when the path of the ``front" point $x$ is an inflectional elastica. We also show the trajectories of points $y_1$ and $y_2$ and depict the motion of the 2-linkage.

Section \ref{sect:un2sol} deals with the case of unequal lengths; in this case, all three of the paths swept out by the three vertices $x, y_1, y_2$ are 2-solitons. Their curvatures satisfy $a E_0 + b E_2 + E_4=0$, with the same parameters $a,b$ for all three points (this does not imply that the three paths are congruent, although that can happen).  The constants $a$ and $b$  are expressible in terms of the phase space variables.  Both are constants of motion for our Hamiltonian flow on the linkage space, with $b$ functionally dependent upon $a$. 

 In Section \ref{sect:int}, we discuss complete integrability of our Hamiltonian system in the equal lengths case.  In Section \ref{sect:fil}, we 
 give a geometric interpretation of the flow of the ``fourth" integral (the above mentioned expression $b$),  a Hamiltonian in its own right, and give a geometric interpretation of its flow in terms of the planar filament equation. 
 
We conclude with Section \ref{sect:concl} that presents a variety of open problems extending, and motivated by, the present research. We illustrate the complexity of some of these problems with computer graphics.

{
This investigation of the geometry of linkages  has involved numerical, graphical, and computer algebra experimentation.  For example,
Proposition \ref{prop:sing2}
 and Theorem \ref{thm:equal} were both conjectured after graphing solutions to the relevant differential equations, followed by rigorous 
proof. 

 In parts of the paper, we lean heavily on computer algebra tools, rather than explicit ``hand'' calculations, for proofs.
Yet, we want the reader to be able to check the details of those calculations.  That is why we provide enough information, so that someone with moderate facility with any one of the standard computer algebra systems can reproduce our results. We look forward to the discovery of proofs  which do not require  the ``brute force" calculations implicit in these computer algebra computations.}

\paragraph{Acknowledgments.} We thank G. Bor, C. Jackman, J. Langer, and M. Levi for interesting discussions and help. RP was supported by the Penn State Shapiro Fund; he is grateful to the Penn State Department of Mathematics for its hospitality. ST was supported by NSF grant DMS-2005444.

\section{Review of 1-linkage geometry} \label{sect:1link}

In this section we briefly review paper \cite{2D} that, along with \cite{BJT}, has motivated the present study.

Consider the space of oriented segments of a fixed (say, unit) length in the plane. The segment can move in such a way that the velocity of its rear end is aligned with the segment. As we mentioned earlier, this is a model of the bicycle kinematics: the rear wheel of the bicycle is fixed on its frame.

The configuration space ${\mathcal C}$ of the segments is 3-dimensional, and the above described non-holonomic constraint defines a 2-dimensional distribution ${\mathcal D}$ therein. This distribution is completely non-integrable, and ${\mathcal C}$ is the space of oriented contact elements in $\R^2$ with its standard contact structure. Due to the Chow-Rashevskii theorem, every pair of points can be connected by a horizontal curve, that is, a curve tangent to  ${\mathcal D}$ (see, e.g., \cite{M}).

The kernel of the differential of the projection ${\mathcal C} \to \R^2$ on the front end of the segment is transverse to ${\mathcal D}$, and this differential  maps the planes of the distribution isomorphically onto $\R^2$.  One gives ${\mathcal D}$ a Riemannian metric by pulling back the Euclidean metric from the plane. The problem is to describe the geodesics of this sub-Riemannian metric. 
We call such geodesics minimizing bicycle paths or bicycling geodesics. Let us reiterate: the length of the bicycle path by definition is the length of the front track.

One uses Hamiltonian formalism to study this problem. Choose an orthonormal frame $(v_1,v_2)$ of ${\mathcal D}$, and consider the vector fields $v_{1}$ and $v_2$ as linear functions $L_{1}$ and $L_2$ on the cotangent bundle $T^*{\mathcal C}$. One obtains a Hamiltonian $H=\frac{1}{2}(L_1^2+L_2^2)$, and the projections to ${\mathcal C}$ of the trajectories of the respective Hamiltonian vector field on the fixed energy surface $H=\frac{1}{2}$ are the arc length parameterized sub-Riemannian geodesics. See \cite{M} for details of this theory.

The main results of \cite{2D} are as follows. 

\begin{thm} \label{thm:5auth}
(i) The front track of a minimizing bicycle path is a straight line or an arc of a non-inflectional elastic curve. Every  shapes of non-inflectional elasticae arise in this way. \\
(ii) An infinitely long bicycle path is a global minimizer (all its segments minimize length between their end points) if either the front track is a straight line and the rear track is either a straight line or a tractrix, or the front track is an Euler soliton and the rear track is a tractrix.
\end{thm}

\begin{figure}[ht]
\centering
\includegraphics[width=1\textwidth]{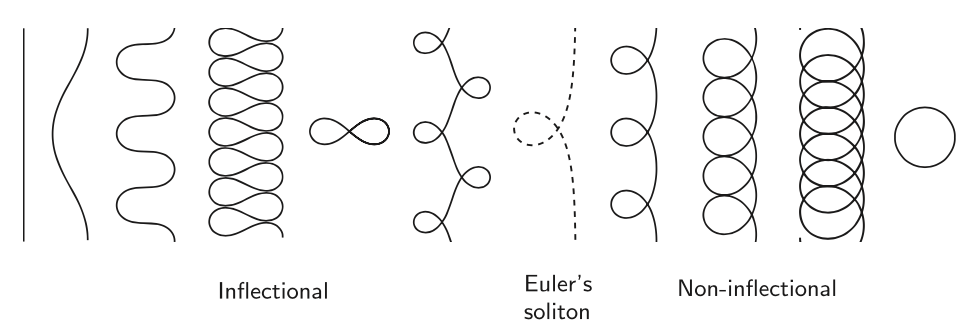}
\caption{A variety of elastic curves.}	
\label{elastica}
\end{figure}

We recall that elastica are the curves that have critical bending energy among the curves of fixed length connecting two fixed points; see, e.g., \cite{Si}. See Figure \ref{elastica}, borrowed from \cite{2D}, for various shapes of elastic curves.

In the bicycle terminology, a tractrix is the rear track when the front track is a straight line.  The Euler soliton is  obtained from the straight line by what is known as the bicycle (or B\"acklund) transformation, consisting of rotation of the bicycle through $180^{\circ}$ about its rear end; this transformation is an isometry of ${\mathcal C}$.

\section{Non-integrability of the distribution} \label{sect:distr}

We start with a description of the configuration space of the linkage and of the distribution defined by the no slip condition.

Introduce the following coordinates: $x=(x_1,x_2)$, the Cartesian coordinates of point $x$, and $\alpha_1, \alpha_2$, the directions of the oriented segments $m_1 x$ and $m_2 x$, respectively. The lengths of the segments $m_1 x$ and $m_2 x$ are denoted by $\ell_1$ and $\ell_2$. The configuration space of the 2-linkage is $\R^2 \times T^2$; let ${\mathcal D}$ be the 2-dimensional distribution defined by the non-holonomic constraint. 

We will further assume that $m_1 \neq m_2$, that is, if $\ell_1=\ell_2$, then $\alpha_1 \neq \alpha_2$. Denote this reduced configuration space by ${\mathcal C}$. As in the 1-linkage case, the kernel of the differential of the projection  ${\mathcal C} \to \R^2$ is transverse to ${\mathcal D}$, and we give the distribution the pullback metric. 

\begin{thm} \label{thm:nonint}
The distribution ${\mathcal D}$ in ${\mathcal C}$ is completely non-integrable with the growth vector $(2,3,4)$.
\end{thm}

\bp
One has
$$
m_1=(x_1-\ell_1\cos\alpha_1,x_2-\ell_1\sin\alpha_1), m_2=(x_2-\ell_2\cos\alpha_2,x_2-\ell_2\sin\alpha_2),
$$
the distribution is given by the Pfaffian equations
$$
\lambda_1=\sin\alpha_1\ dx_1 - \cos\alpha_1\ dx_2 + \ell_1 d\alpha_1=0,\ \lambda_2=\sin\alpha_2\ dx_1 - \cos\alpha_2\ dx_2 + \ell_2 d\alpha_2=0, 
$$
and a basis of horizontal fields is given by
\begin{equation} \label{eq:base}
v_1=\partial_{x_1}-\frac{\sin\alpha_1}{\ell_1}\ \partial_{\alpha_1}-\frac{\sin\alpha_2}{\ell_2}\ \partial_{\alpha_2},\ 
v_2=\partial_{x_2}+\frac{\cos\alpha_1}{\ell_1}\ \partial_{\alpha_1}+\frac{\cos\alpha_2}{\ell_2}\ \partial_{\alpha_2}.
\end{equation}
Then
\begin{equation*}
\begin{split}
&[v_1,v_2]=\frac{1}{\ell_1^2}\ \partial_{\alpha_1}+\frac{1}{\ell_2^2}\ \partial_{\alpha_2}, 
[[v_1,v_2],v_1]=-\frac{\cos\alpha_1}{\ell_1^3}\ \partial_{\alpha_1} - \frac{\cos\alpha_2}{\ell_2^3}\ \partial_{\alpha_2},\\
&[[v_1,v_2],v_2]=-\frac{\sin\alpha_1}{\ell_1^3}\ \partial_{\alpha_1} - \frac{\sin\alpha_2}{\ell_2^3}\ \partial_{\alpha_2}.
\end{split}
\end{equation*}
If $\ell_1\neq \ell_2$, then the matrix
\begin{equation} \label{eq:mat}
\begin{pmatrix}
\frac{1}{\ell_1^2} & \frac{\cos\alpha_1}{\ell_1^3} & \frac{\sin\alpha_1}{\ell_1^3}\\
\frac{1}{\ell_2^2} & \frac{\cos\alpha_2}{\ell_2^3} & \frac{\sin\alpha_2}{\ell_2^3}
\end{pmatrix}
\end{equation}
has rank 2. Indeed, if
$$
\left(\frac{\cos\alpha_1}{\ell_1^3}, \frac{\cos\alpha_2}{\ell_2^3}\right) = t \left(\frac{1}{\ell_1^2}, \frac{1}{\ell_2^2} \right),
\left(\frac{\sin\alpha_1}{\ell_1^3}, \frac{\sin\alpha_2}{\ell_2^3}\right) = s \left(\frac{1}{\ell_1^2}, \frac{1}{\ell_2^2} \right),
$$
then $(t^2+s^2)\ell_1^2=(t^2+s^2)\ell_2^2$, and hence $\ell_1=\ell_2$.

Since the matrix has rank 2, both $\partial_{\alpha_1}$ and $\partial_{\alpha_2}$ are linear combinations of $[v_1,v_2], [[v_1,v_2],v_1]$ and $[[v_1,v_2],v_2]$, and then both $\partial_{x_1}$ and $\partial_{x_2}$ are linear combinations of these vectors, and $v_1$ and $v_2$. 

If $\ell_1=\ell_2$, we may assume that they are equal to 1. In this case the rank of matrix (\ref{eq:mat}) is not equal to 2 only if $\alpha_1=\alpha_2$, which is excluded.  In addition, 
$$
v_2+[[v_1,v_2],v_1]=\partial_{x_2},\ v_1-[[v_1,v_2],v_2]=\partial_{x_2},
$$
as needed.
\ep

The Chow-Rashevskii theorem implies that any pair of points of ${\mathcal C}$ can be connected by a horizontal curve. 

A 2-dimensional distribution on a 4-dimensional manifold with the growth vector $(2,3,4)$ is called an Engel structure. Like contact structures, all Engel structures are locally diffeomorphic. The normal form of the distribution is given, in local coordinates $(x,y,z,w)$, by the two 1-forms $dz-ydx$ and $dy-wdx$, see \cite{M}.

\section{Singular curves} \label{sect:sing}

A specific phenomenon of sub-Riemannian geometry is that the space of horizontal paths that connect two fixed points may have singularities. These singular curves can be sub-Riemannian geodesics that are not detected by the Hamiltonian formalism of
Section \ref{sect:eq};
see \cite{M,M2}.  

It is known that Englel manifolds are foliated by singular curves: in the above mentioned normal form, they are the trajectories of the vector field $\partial/\partial w$. Let us describe the singular curves in our setting.

\begin{thm} \label{thm:sing}
The singular curves in ${\mathcal C}$ are the integral curves of the vector field
$$
\xi=\left(\frac{\sin\alpha_1}{\ell_1}-\frac{\sin\alpha_2}{\ell_2}\right) v_1 - \left(\frac{\cos\alpha_1}{\ell_1}-\frac{\cos\alpha_2}{\ell_2}\right) v_2,
$$
where $v_1$ and $v_2$ are as in (\ref{eq:base}).
\end{thm}

\bp
Set ${\mathcal D}^2={\mathcal D}+[{\mathcal D},{\mathcal D}]$. 
A desired vector field $\xi$ is characterized by the property that $[\xi, {\mathcal D}^2] \subset {\mathcal D}^2$, see \cite{M2}. 

The differential form $\lambda:=\ell_1 \lambda_1-\ell_2 \lambda_2$ annihilates ${\mathcal D}^2$.
Let $\xi=f v_1 + g v_2$. Then the fields $[\xi,v_1]$ and $[\xi,v_2]$ are already in ker $\lambda$. One has
$$
\lambda\left(\left[\frac{1}{\ell_1^2}\ \partial_{\alpha_1}+\frac{1}{\ell_2^2}\ \partial_{\alpha_2},\xi     \right]\right)=
f \left(\frac{\cos\alpha_1}{\ell_1} - \frac{\cos\alpha_2}{\ell_2}  \right) + g \left(\frac{\sin\alpha_1}{\ell_1} - \frac{\sin\alpha_2}{\ell_2}  \right).
$$
This must vanish, hence one can choose
$$
f= \frac{\sin\alpha_1}{\ell_1}-\frac{\sin\alpha_2}{\ell_2},\ g=- \frac{\cos\alpha_1}{\ell_1}+\frac{\cos\alpha_2}{\ell_2},
$$
as claimed.
\ep

The explicit formula is as follows:
\begin{equation*}
\begin{split}
\xi&=\left(\frac{\sin\alpha_1}{\ell_1} - \frac{\sin\alpha_2}{\ell_2}  \right)\ \partial_{x_1} - \left(\frac{\cos\alpha_1}{\ell_1} - \frac{\cos\alpha_2}{\ell_2}  \right)\ \partial_{x_2}  \\
 &-\left( \frac{1}{\ell_1^2} - \frac{\cos(\alpha_1-\alpha_2)}{\ell_1\ell_2} \right)\ \partial_{\alpha_1} +  
\left( \frac{1}{\ell_2^2} - \frac{\cos(\alpha_1-\alpha_2)}{\ell_1\ell_2} \right)\ \partial_{\alpha_2}.
\end{split}
\end{equation*}

\begin{prop} \label{prop:sing}
If $\ell_1=\ell_2$, then the motion along the singular curves is as follows: the angle $(\alpha_1+\alpha_2)/2$ remains constant, and point $x$ moves along a straight line that bisects the angle between the segments $xm_1$ and $xm_2$. The points $m_1$ and $m_2$ move along the tractrices which are symmetric with respect to this line. This singular curve is a minimizing geodesic.
\end{prop}

\bp
If $\ell_1=\ell_2=1$, the formula simplifies to
\begin{align*}
&\xi  = 2 \sin\left(\frac{\alpha_1-\alpha_2}{2} \right) \times\\
&\left[ \cos\left(\frac{\alpha_1  + \alpha_2}{2}\right) \partial_{x_1} + \sin\left(\frac{\alpha_1+\alpha_2}{2}\right) \partial_{x_2}+ \sin\left(\frac{\alpha_1-\alpha_2}{2}\right) (\partial_{\alpha_2} - \partial_{\alpha_1})  \right].
\end{align*}
The factor in front does not vanish on ${\mathcal C}$, and we arrive at the differential equations describing the singular curves:
\begin{align*}
&\dot x_1=\cos\left(\frac{\alpha_1+\alpha_2}{2}\right), \dot x_2=\sin\left(\frac{\alpha_1+\alpha_2}{2}\right),\\
&\dot\alpha_1 = - \sin\left(\frac{\alpha_1-\alpha_2}{2}\right) , \dot\alpha_2=\sin\left(\frac{\alpha_1-\alpha_2}{2}\right) .
\end{align*}
It follows that $\dot\alpha_1+\dot\alpha_2=0$, hence $\alpha_1+\alpha_2$ is constant, and then 
$$
\dot x = \left(\cos\left(\frac{\alpha_1+\alpha_2}{2}\right),\sin\left(\frac{\alpha_1+\alpha_2}{2}\right)\right),
$$
as claimed.

The last claim follows from the fact the the trajectory of point $x$ is a straight line.
\ep

We also have the following

\begin{prop} \label{prop:sing2}
If $\ell_1\neq \ell_2$, then the $x$-projection of a singular curve on the plane is an elastic curve.
\end{prop}

\begin{figure}[ht]
\centering
\includegraphics[width=.45\textwidth]{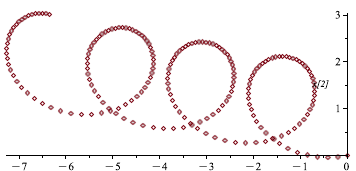}
\includegraphics[width=.5\textwidth]{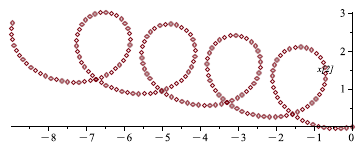}
\caption{Left: planar projection of a singular geodesic; right: the same curve parameterized by arc length.}	
\label{para}
\end{figure}

\bp
We reiterate the differential equations for the singular geodesic:
\begin{align*}
\frac{d}{d t}\alpha_{1}\! \left(t \right) &= 
\frac{{\ell}_{1} \cos \! \left(-\alpha_{2}\! \left(t \right)+\alpha_{1}\! \left(t \right)\right)-{\ell}_{2}}{{\ell}_{1}^{2} {\ell}_{2}} \\
\frac{d}{d t}\alpha_{2}\! \left(t \right)  &= 
-\frac{{\ell}_{2} \cos \! \left(-\alpha_{2}\! \left(t \right)+\alpha_{1}\! \left(t \right)\right)-{\ell}_{1}}{{\ell}_{1} {\ell}_{2}^{2}}\\
\frac{d}{d t}x_{1}\! \left(t \right)  &= 
\frac{\sin \! \left(\alpha_{1}\! \left(t \right)\right) {\ell}_{2}-\sin \! \left(\alpha_{2}\! \left(t \right)\right) 
{\ell}_{1}}{{\ell}_{1} {\ell}_{2}}\\
\frac{d}{d t}x_{2}\! \left(t \right) &= 
-\frac{\cos \! \left(\alpha_{1}\! \left(t \right)\right) \mathit{\ell}_{2}-\cos \! \left(\alpha_{2}\! \left(t \right)\right) \mathit{\ell}_{1}}{\mathit{\ell}_{1} \mathit{\ell}_{2}}
\end{align*}
We experiment, and plot the planar projection of the solution curve  (the following  observations  are independent of the particular choice of initial conditions), see Figure \ref{para}.

The curve looks like an elastica and, by plotting points equally spaced in the parameter $t$, we see that it is not unit speed
parameterized. 

  We will fix this, and then compute the curvature of the planar curve. Let $S$ be the speed; it is easily computed to be:
$$
S = \sqrt{-\frac{2 {\ell}_{1} {\ell}_{2} \cos \! \left(-\alpha_{2}\! \left(t \right)+\alpha_{1}\! \left(t \right)\right)-{\ell}_{1}^{2}-{\ell}_{2}^{2}}{{\ell}_{1}^{2} {\ell}_{2}^{2}}}.
$$
Dividing our set of equations by $S$, we obtain a time-scaled version of the differential equations, with the same solution curves (reparameterized). We plot solutions to these modified equations, and obtain the same planar curve, but unit speed parameterized,
see Figure \ref{para}.

Let ${\bf p} = (x_1(t), x_2(t))$, the curve defined by the {second} pair of differential equations.  By  scaling,
we have $\dot p = T$, the unit tangent to the curve.  One immediately obtains the unit normal $N$, and the computation
$< \dot T,N>$  results in the curvature.  We differentiate the curvature twice:
\begin{align*}
\kappa  &= 
\frac{\mathit{\ell}_{1}^{2}-\mathit{\ell}_{2}^{2}}{\mathit{\ell}_{1} \mathit{\ell}_{2} \sqrt{-2 \mathit{\ell}_{1} \mathit{\ell}_{2} \cos \! \left(-\alpha_{2}\! +\alpha_{1}\! \right)+\mathit{\ell}_{1}^{2}+\mathit{\ell}_{2}^{2}}}, \\
\dot\kappa   &=
\frac{\left(-\mathit{\ell}_{1}^{2}+\mathit{\ell}_{2}^{2}\right) \sin \! \left(-\alpha_{2}\! +\alpha_{1}\! \right)}{\mathit{\ell}_{1} \mathit{\ell}_{2} \left(2 \mathit{\ell}_{1} \mathit{\ell}_{2} \cos \! \left(-\alpha_{2}\! +\alpha_{1}\! \right)-\mathit{\ell}_{1}^{2}-\mathit{\ell}_{2}^{2}\right)}.\\
\ddot\kappa &= 
\frac{\cos \! \left(-\alpha_{2}\! +\alpha_{1}\! \right) \mathit{\ell}_{1}^{4}-\cos \! \left(-\alpha_{2}\! +\alpha_{1}\! \right) \mathit{\ell}_{2}^{4}-2 \mathit{\ell}_{1}^{3} \mathit{\ell}_{2}+2 \mathit{\ell}_{1} \mathit{\ell}_{2}^{3}}{\mathit{\ell}_{1}^{2} \mathit{\ell}_{2}^{2} \left(2 \mathit{\ell}_{1} \mathit{\ell}_{2} \cos \! \left(-\alpha_{2}\! +\alpha_{1}\! \right)-\mathit{\ell}_{1}^{2}-\mathit{\ell}_{2}^{2}\right) \sqrt{-2 \mathit{\ell}_{1} \mathit{\ell}_{2} \cos \! \left(-\alpha_{2}\! +\alpha_{1}\! \right)+\mathit{\ell}_{1}^{2}+\mathit{\ell}_{2}^{2}}},
\end{align*}
and compute the expression 
$
\kappa'' + \frac{\kappa^3}{2} + A \kappa$. It vanishes if $A=-\frac{\ell_{1}^{2}+\ell_{2}^{2}}{2 \ell_{1}^{2} \ell_{2}^{2}}$, hence the curve is an elastica.
\ep

\section{Geodesic equations for 2-linkages} \label{sect:eq}
Now we apply the Hamiltonian formalism to deduce the equations for non-singular geodesics.

An orthonormal basis of horizontal vector fields is given by $(v_1,v_2)$ from equation (\ref{eq:base}). Consider the following coordinates in $T^*{\mathcal C}$:
$$
(x_1, x_2, \alpha_1, \alpha_2, p_1, p_2, \eta_1, \eta_2),
$$
where the last four are the fiber coordinates conjugated to the first four, respectively. Then the vector fields $v_1$ and $v_2$ become linear functions on $T^*{\mathcal C}$:
\begin{equation} \label{eq:fnct}
L_1=p_1-\frac{1}{\ell_1}\eta_1 \sin\alpha_1 -\frac{1}{\ell_2}\eta_2 \sin\alpha_2,\ 
L_2=p_2+\frac{1}{\ell_1}\eta_1 \cos\alpha_1 +\frac{1}{\ell_2}\eta_2 \cos\alpha_2.
\end{equation}
The Hamiltonian is given by
$$
H=\frac{1}{2} (L_1^2+L_2^2).
$$
Set $H=\frac{1}{2}$, and introduce the angle $\gamma$  by setting $L_1=\cos\gamma, L_2=\sin\gamma$.

The Hamiltonian equations read (after an application of appropriate trigonometric identities):
\begin{equation} \label{eq:Ham}
\begin{split}
&\dot x_1=\cos\gamma, \dot x_2 = \sin\gamma, \dot \alpha_1 = \frac{\sin(\gamma-\alpha_1)}{\ell_1}, \dot \alpha_2 = \frac{\sin(\gamma-\alpha_2)}{\ell_2},\\
&\dot p_1 = \dot p_2 =0, 
\dot \eta_1 = \frac{\eta_1 \cos(\gamma-\alpha_1)}{\ell_1},  \dot \eta_2 = \frac{\eta_2 \cos(\gamma-\alpha_2)}{\ell_2},
\end{split}
\end{equation}
where dot is the time derivative.

We are interested in the curve $x(t)$. Note that $t$ is the arc length parameter;  let $\kappa(t)$ be the curvature of this curve. 

\begin{prop} \label{prop:curv}
One has
\begin{equation} \label{eq:curv1}
\begin{split}
\kappa=\frac{\eta_1}{\ell_1^2}+\frac{\eta_2}{\ell_2^2},\ \dot \kappa = \frac{\eta_1 \cos(\gamma-\alpha_1)}{\ell_1^3} + \frac{\eta_2 \cos(\gamma-\alpha_2)}{\ell_2^3},\\ 
\ddot \kappa = \frac{\eta_1}{\ell_1^4}+\frac{\eta_2}{\ell_2^4}- \kappa \left[\frac{\eta_1\sin(\gamma-\alpha_1)}{\ell_1^3}+\frac{\eta_2\sin(\gamma-\alpha_2)}{\ell_2^3}\right].
\end{split}
\end{equation}
\end{prop}

\bp
One has $\kappa=\dot\gamma$. 
Differentiate (say) the second equation in (\ref{eq:fnct}), using that $\dot p_2=0$. Use equations (\ref{eq:Ham}) and some trigonometry to obtain 
$$
\dot \gamma \sin\gamma= \frac{\eta_1}{\ell_1^2} \sin\gamma + \frac{\eta_2}{\ell_2^2} \sin\gamma,
$$ 
which implies the equation for the curvature. The other two equations are obtained by differentiating the first one, substituting the values of the derivatives from (\ref{eq:Ham}), and using suitable  trigonometric identities.
\ep

\begin{rmrk} \label{rmk:bic}
The equations 
$$
\dot \alpha_1 = \frac{\sin(\gamma-\alpha_1)}{\ell_1}, \dot \alpha_2 = \frac{\sin(\gamma-\alpha_2)}{\ell_2}
$$
in (\ref{eq:Ham}) are just the ``bicycle equations" expressing the no-skid non-holonomic constraints on the points $m_1$ and $m_2$, see, e.g., \cite{BLPT,FLT}.
\end{rmrk}

\paragraph{Lifting 1-linkage geodesics to 2-linkage geodesics} 
Consider 1-linkage $x m_1$, and let $x(t) m_1(t), t\in[0,T]$, be a bicycling geodesic. Let $m_2=m_2(0)$ be another point, interpreted as the rear wheel of the bicycle $x m_2$. The no-skid  constraint uniquely determines the motion $m_2(t) , t\in[0,T]$, and we obtain a horizontal path in the 2-linkage configuration space ${\mathcal C}$.

\begin{prop} \label{prop:lift}
This path is geodesic. 
\end{prop} 

\bp
We need to show that the path $x(t) m_1(t) m_2(t)$ minimizes the distance between its two sufficiently close points. 

Assume that there exists a shorter horizontal path from $x(t_1) m_1(t_1) m_2(t_1)$ to $x(t_2) m_1(t_2) m_2(t_2)$. Forgetting point $m_2$, we obtain a shorter horizontal path in the space of 1-linkage connecting $x(t_1) m_1(t_1)$ and $x(t_2) m_1(t_2)$. This contradicts the assumption that $x(t) m_1(t)$ is a bicycling geodesic.
\ep

If $x(t) m_1(t)$ is a bicycling geodesic, then the curve $x(t)$ is a non-inflectional elastica, see Section \ref{sect:1link}. The above proposition implies that all non-inflectional elasticae  are  the projections of some geodesics in ${\mathcal C}$ to $\R^2$. It also follows that there are infinitely many lifts of a non-inflectional elastica to a geodesic curve in  ${\mathcal C}$.

\section{The case of equal  length line segments and the appearance of elastica} \label{sect:ela}
In this section we assume that $\ell_1=\ell_2$ and, by scaling, that this length is unit. 

\begin{thm} \label{thm:equal}
The projection $x(t)$ of a geodesic in ${\mathcal C}$ to $\R^2$ is an elastic curve. 
\end{thm}

\bp
Elasticae \ are characterized by the differential equation on their curvature as a function of the arc length parameter
$$
\ddot \kappa + \frac{1}{2} \kappa^3 + A \kappa =0,
$$
so (ignoring the case of a straight line where the curvature is identically zero) we need to show that 
$$
\frac{\ddot \kappa}{\kappa} + \frac{1}{2} \kappa^2 = const.
$$

Equation (\ref{eq:curv1}) implies that 
$$
\frac{\ddot \kappa}{\kappa} = 1 - {\eta_1\sin(\gamma-\alpha_1)}-{\eta_2\sin(\gamma-\eta_2)}.
$$
One can solve equations (\ref{eq:fnct}) for $\eta_1$ and $\eta_2$
\begin{equation} \label{eq:elim}
\eta_1=\frac{\cos(\gamma-\alpha_2)-p_1\cos\alpha_2-p_2\sin\alpha_2}{\sin(\alpha_2-\alpha_1)},\ 
\eta_2=\frac{\cos(\gamma-\alpha_1) -p_1\cos\alpha_1-p_2\sin\alpha_1}{\sin(\alpha_1-\alpha_2)}
\end{equation}
and substitute to the above equation to obtain, using trigonometric identities:
\begin{equation}  \label{eq:doubleprime}
\frac{\ddot \kappa}{\kappa} = p_1 \cos\gamma+p_2\sin\gamma.
\end{equation}

We want to check that
$$
\frac{d}{dt}  \left(\frac{\ddot \kappa}{\kappa}\right) + \kappa\dot \kappa =0.
$$
One has
$$
\frac{d}{dt} (p_1 \cos\gamma+p_2\sin\gamma) = (-p_1\sin\gamma+p_2\cos\gamma) \kappa,
$$
so it remains to check that 
\begin{equation} \label{eq:prime}
-p_1\sin\gamma+p_2\cos\gamma = -\dot \kappa.
\end{equation}
To see this, take the value of $\dot \kappa$ from (\ref{eq:curv1}) and substitute $\eta_1$ and $\eta_2$ from (\ref{eq:elim}) to obtain a true identity. This shows that $x(t)$ is an elastic curve.
\ep


Let us write the equations of elastica in two forms:
\begin{equation} \label{eq:forms}
\frac{\ddot \kappa}{\kappa} + \frac{1}{2} \kappa^2 + A=0,\ \dot \kappa^2+\frac{1}{4} \kappa^4 + A\kappa^2+B=0,
\end{equation}
the second one is obtained from the first by integration.

\begin{lemma} \label{lm:AB}
One has $A^2-B= p_1^2+p_2^2.$
\end{lemma} 

\bp
Equations (\ref{eq:doubleprime}) and (\ref{eq:prime}) imply
\begin{equation*} \label{eq:rel}
\left(\frac{\ddot \kappa}{\kappa}\right)^2+\dot \kappa^2=p_1^2+p_2^2.
\end{equation*}
Hence
$$
\left( \frac{1}{2} \kappa^2 + A\right)^2- \frac{1}{4} \kappa^4 - A\kappa^2-B = p_1^2+p_2^2,
$$
therefore
$$
A^2-B= p_1^2+p_2^2,
$$
as claimed.
\ep

Set 
$$
G= -A = \frac{\ddot \kappa}{\kappa} + \frac{1}{2} \kappa^2.
$$
As we have seen, this is another integral of motion, in addition to $H, p_1$, and $p_2$. We shall address the independence of these integrals in Section \ref{sect:int}. 

Using (\ref{eq:doubleprime}), we also obtain a differential equation on $\gamma$:
$$
\dot \gamma = \kappa= \sqrt{2(G - p_1\cos\gamma-p_2\sin\gamma)}.
$$
The substitution $y=\tan(\gamma /2)$ transforms this to the equation $\dot y = \sqrt{Q(y)}$,
where $Q$ is a polynomial of degree 4.

Note that all non-closed elastic curves are contained in a strip, see Figure \ref{elastica}. Let us call the direction of the strip the direction of the respective elastica. 

\begin{lemma} \label{lm:dir}
The vector $(p_1,p_2)$ is parallel to the direction of the respective elastica.
\end{lemma}

\bp
Write $(p_1,p_2)=r (\cos \phi,\sin\phi)$, so $\phi$ is the direction of this vector. Then equation (\ref{eq:prime}) can be written as $\dot \kappa = r \sin(\gamma-\phi)$. Hence if $\dot \kappa =0$, then $\phi=\gamma \mod \pi$. That is, $\phi$ is the direction of the elastica at its vertices (the point where the curvature is extremal). This direction  coincides with the direction of the elastica, see \cite{Si} and Figure \ref{elastica}. 
\ep

This lemma agrees with the fact that when the elastica is a circle (all points are vertices), one has $p_1=p_2=0$. We also note that 
$r=|(p_1,p_2)|$  is determined by the function $\g(t)$: 
$$
r=\left|\frac{\dot \kappa}{\sin(\gamma-\phi)}\right|=\left|\frac{\ddot \g}{\sin(\gamma-\phi)}\right|.
$$

One can reduce the number of variables in equations (\ref{eq:Ham}) by solving equations (\ref{eq:fnct}) for $\eta_1$ and $\eta_2$, and substituting the result in (\ref{eq:Ham}). This yields the system of equations
\begin{equation} \label{eq:abg}
\begin{aligned}
&\dot \alpha_1 = {\sin(\gamma-\alpha_1)},\ \dot \alpha_2 = {\sin(\gamma-\alpha_2)}, \\
&\dot \gamma =  \frac{\sin\left(\gamma-\frac{\alpha_1+\alpha_2}{2} \right) +p_1 \sin\left(\frac{\alpha_1+\alpha_2}{2} \right) - p_2 \cos\left(\frac{\alpha_1+\alpha_2}{2} \right)}{\cos\left(\frac{\alpha_2-\alpha_1}{2} \right)}.
\end{aligned}
\end{equation}

As we saw in Section \ref{sect:eq}, and in view of Section \ref{sect:1link}, every non-inflectional elastica is a projection of a geodesic curve. What about inflectional elasticae?

\begin{thm} \label{thm:all}
Every shape of elastica appears as the projection of a geodesic curve in ${\mathcal C}$.
\end{thm}

\begin{proof}
Given  an arc length  parameterized elastica, we have an explicit function $\gamma(t)$ (given in terms of elliptic functions). Then the constants $p_1$ and $p_2$ are determined -- see Lemma \ref{lm:dir} and the discussion after it.   

Consider equations (\ref{eq:abg}). Once the initial values $\alpha_1(0)$ and $\alpha_2(0)$ are chosen, the first two differential equations determine the functions $\alpha_1(t)$ and $\alpha_2(t)$ uniquely. 

Suppose that $\alpha_1(0)$ and $\alpha_2(0)$ are chosen in such a way that the third equation (\ref{eq:abg}) holds for $t=0$. We claim that then it holds for all values of $t$. This would imply that the elastica under consideration is the projection of the respective geodesic.

Set
$$
f(t) = \dot \gamma -  \frac{\sin\left(\gamma-\frac{\alpha_1+\alpha_2}{2} \right) +p_1 \sin\left(\frac{\alpha_1+\alpha_2}{2} \right) - p_2 \cos\left(\frac{\alpha_1+\alpha_2}{2} \right)}{\cos\left(\frac{\alpha_2-\alpha_1}{2} \right)}.
$$
Using the differential equations for $\alpha_1(t)$ and $\alpha_2(t)$, and equation (\ref{eq:prime}), one calculates that $f$
satisfies the following first order linear differential equation
$$
\dot f = - \frac{\cos\left(\gamma-\frac{\alpha_1+\alpha_2}{2}\right) }{ \cos\left(\frac{\alpha_1-\alpha_2}{2}\right)}  f. 
$$
Since $f(0)=0$ by assumption, $f(t)\equiv 0$ is a solution. Then the standard theorem of existence and uniqueness of solutions to ordinary differential equations (e.g.,\cite{CL}) implies that the third equation (\ref{eq:abg}) holds for all $t$, as needed.
\end{proof}

The shape of the elastica satisfying equation (\ref{eq:forms}) is determined by the scale-invariant parameter $\mu=B/A^2$ (our choice of constant $B$  is slightly different from that in \cite{2D}). One has
$$
\mu=1 - \frac{p_1^2+p_2^2}{G} \le 1.
$$
A non-inflectional elastica corresponds to $\mu > 0$, and  an inflectional elastica to $\mu<0$. 

Thus, our situation differs from that occuring in \cite{2D};  for the equal lengths case of 2-linkages, $\mu$ may be negative; this  fact allows for inflectional
elasticae.  In Example \ref{ex:inflection}, we show a convenient way to generate all of the inflectional elasticae.

\section{Computer-assisted proofs} \label{sect:comp}

 It is perhaps not surprising that calculations for 2-linkages with unequal length line segments are harder than in the equal length case.  
 When we  have been unable to prove results ``by hand'', we have resorted to computer algebra calculations.
 
In particular, we make elementary but extensive use of a particular tool in computer algebra, the use of {\it Gr\"obner} bases, which allows us to use ``brute force"  to do the relevant calculations.  Gr\"obner basis tools exist for all of the major computer algebra systems (in particular, the popular tools, Maple and Mathematica). The paper \cite{St} gives a nice short introduction, but we say a few words here.

Given an ideal $I$ in a polynomial ring $R = F[x_1, x_2, \dots, x_n] $ with generating polynomials 
$PB =\{p_1, p_2 , \dots , p_k\}$, how can one 
decide in an algorithmic manner whether another polynomial $p$ is an element of the ideal?  And, if $p$ is not an element
of {\it I}, is there a canonical representative of $p$ in the quotient ring $R \mod I$?  This question was addressed by Buchberger; the algorithm requires the computation of an alternative set of generators for $I$  (the {\it Gr\"obner basis} $PGB$), and then a computation using this new basis to find the canonical representative of $p \mod I$. See, e.g., \cite{Laza}. 

Of course, in order for this to be useful, one must phrase the problem of interest entirely in terms  of polynomial calculations.   We show how this applies to a second proof that, for equal lengths, the curve $x(t)$ is an elastic curve. We describe this computation in excruciating detail, and we will be more terse later in the paper.

The idea is to express all quantities and equations as polynomial expressions in the  set of variables
$$vars =\{ \alpha_1,\alpha_2, \eta_1, \eta_2, p_1, p_2, \cca, \ccb, \ssa, \ssb\}, $$
where $\cca, \ccb, \ssa, \ssb$ are polynomial variables corresponding to the functions
$\cos\alpha_1, \, \cos\alpha_2, \, \sin\alpha_1, \, \sin\alpha_2$ respectively.  For example,
the polynomial representing the Hamiltonian can be written as
\begin{align*}
H  & \longleftrightarrow
Hpoly =-\frac{\mathit{\ssa} \eta_{1} p_{1}}{\mathit{\ell}_{1}}-\frac{\mathit{\ssb} \eta_{2} p_{1}}{\mathit{\ell}_{2}}+\frac{\eta_{1} \mathit{\cca} p_{2}}{\mathit{\ell}_{1}}+\frac{\eta_{2} \mathit{\ccb} p_{2}}{\mathit{\ell}_{2}}+\frac{p_{1}^{2}}{2}+\frac{p_{2}^{2}}{2} \\
&+\frac{\eta_{1} \mathit{\cca} \eta_{2} \mathit{\ccb}}{\mathit{\ell}_{1} \mathit{\ell}_{2}}+\frac{\mathit{\ssa} \eta_{1} \mathit{\ssb} \eta_{2}}{\mathit{\ell}_{1} \mathit{\ell}_{2}}+\frac{\eta_{2}^{2}}{2 \mathit{\ell}_{2}^{2}}+\frac{\eta_{1}^{2}}{2 \mathit{\ell}_{1}^{2}}.
\end{align*}
We emphasize that  $H$ depends on the functions $\alpha_1(t), \alpha_2(t), \dots$, whereas the expression on the right is a polynomial in the variables 
$vars$.  There is a correspondence, but  they are not equal, hence the $\longleftrightarrow$ symbol.

The geodesic equations 
can be written as follows:
\begin{equation}
\begin{aligned}  \label{eq:bigdes}
\frac{d}{d t}\alpha_{1}\!  &= 
\frac{\cos \! \left(\alpha_{1}\! \right) p_{2}\! }{\mathit{\ell}_{1}}-\frac{p_{1}\!  \sin \! \left(\alpha_{1}\! \right)}{\mathit{\ell}_{1}}\\
& +\frac{\eta_{2}\!  \cos \! \left(\alpha_{1}\! \right) \cos \! \left(\alpha_{2}\! \right)}{\mathit{\ell}_{1} \mathit{\ell}_{2}}+\frac{\eta_{2}\!  \sin \! \left(\alpha_{1}\! \right) \sin \! \left(\alpha_{2}\! \right)}{\mathit{\ell}_{1} \mathit{\ell}_{2}}+\frac{\eta_{1}\! }{\mathit{\ell}_{1}^{2}},\\
\frac{d}{d t}\alpha_{2}\!  & = 
\frac{\cos \! \left(\alpha_{2}\! \right) p_{2}\! }{\mathit{\ell}_{2}}-\frac{p_{1}\!  \sin \! \left(\alpha_{2}\! \right)}{\mathit{\ell}_{2}}\\
& +\frac{\eta_{1}\!  \cos \! \left(\alpha_{1}\! \right) \cos \! \left(\alpha_{2}\! \right)}{\mathit{\ell}_{1} \mathit{\ell}_{2}}+\frac{\eta_{1}\!  \sin \! \left(\alpha_{1}\! \right) \sin \! \left(\alpha_{2}\! \right)}{\mathit{\ell}_{1} \mathit{\ell}_{2}}+\frac{\eta_{2}\! }{\mathit{\ell}_{2}^{2}},\\
\frac{d}{d t}\eta_{1}\!  & = 
\frac{\eta_{1}\!  \sin \! \left(\alpha_{1}\! \right) p_{2}\! }{\mathit{\ell}_{1}}+\frac{p_{1}\!  \cos \! \left(\alpha_{1}\! \right) \eta_{1}\! }{\mathit{\ell}_{1}}\\
&+\frac{\eta_{1}\!  \eta_{2}\!  \sin \! \left(\alpha_{1}\! \right) \cos \! \left(\alpha_{2}\! \right)}{\mathit{\ell}_{1} \mathit{\ell}_{2}}-\frac{\eta_{1}\!  \eta_{2}\!  \cos \! \left(\alpha_{1}\! \right) \sin \! \left(\alpha_{2}\! \right)}{\mathit{\ell}_{1} \mathit{\ell}_{2}},\\
\frac{d}{d t}\eta_{2}\!  &= 
\frac{\eta_{2}\!  \sin \! \left(\alpha_{2}\! \right) p_{2}\! }{\mathit{\ell}_{2}}+\frac{p_{1}\!  \cos \! \left(\alpha_{2}\! \right) \eta_{2}\! }{\mathit{\ell}_{2}}\\
&-\frac{\eta_{1}\!  \eta_{2}\!  \sin \! \left(\alpha_{1}\! \right) \cos \! \left(\alpha_{2}\! \right)}{\mathit{\ell}_{1} \mathit{\ell}_{2}}
 +\frac{\eta_{1}\!  \eta_{2}\!  \cos \! \left(\alpha_{1}\! \right) \sin \! \left(\alpha_{2}\! \right)}{\mathit{\ell}_{1} \mathit{\ell}_{2}},\\
\frac{d}{d t}p_{1}\!  & = 0, \quad
\frac{d}{d t}p_{2}\!  = 0, \\
\frac{d}{d t}x_{1}\!  & = 
-\frac{\sin \! \left(\alpha_{1}\! \right) \eta_{1}\! }{\mathit{\ell}_{1}}-\frac{\sin \! \left(\alpha_{2}\! \right) \eta_{2}\! }{\mathit{\ell}_{2}}+p_{1}\! ,\\
\frac{d}{d t}x_{2}\!  & = 
\frac{\eta_{1}\!  \cos \! \left(\alpha_{1}\! \right)}{\mathit{\ell}_{1}}+\frac{\eta_{2}\!  \cos \! \left(\alpha_{2}\! \right)}{\mathit{\ell}_{2}}+p_{2}\! .  \\
\end{aligned}
\end{equation}

We make note  of a trivial, but important technical point that could be troublesome.  The relevant computer tools
(for example, {\it Basis, NormalForm} in Maple, {\it GroebnerBasis, PolynomialReduce} in Mathematica)
 are used to manipulate {\it polynomials } in 
the {\it variables} 
$$vars =\{ \alpha_1,\alpha_2, \eta_1, \eta_2, p_1, p_2, \cca, \ccb, \ssa, \ssb\}, $$
 whereas the differentiations will be applied to {\it polynomials}
 in the {\it functions}
$$\{x_1(t),  \,   x_2(t),  p_1(t), p _2(t),  \cos(\alpha_1(t)), \, \cos(\alpha_2(t)), \, \sin(\alpha_1(t)), 
\, \sin(\alpha_2(t)), \eta_1(t), \eta_2(t)\}.$$
It is trivial  to create tools {\it FtoP, PtoF}  which translate one to the other.

We finally begin the computation.  We will be interested in working with the condition $H - \frac{1}{2}=0 $, so that the 
projected curve $x = (x_1(t), x_2(t))$ is arc length parameterized.  Thus, our ideal $I $ will be generated by $Hpoly -\frac{1}{2}$, 
as well as the trivial trigonometric identities:
$$PB = \{\mathit{Hpoly}-\frac{1}{2}, \mathit{\cca}^{2}+\mathit{\ssa}^{2}-1, 
\mathit{\ccb}^{2}+\mathit{\ssb}^{2}-1\}.$$
Using the appropriate tool, one computes its associated Gr\"obner basis $PGB$, allowing us to do polynomial calculations $\mod I$.
We re-derive the formula for curvature: we have
\begin{align*}
\frac{d}{d t}x_{1}\! \left(t \right) & \leftrightarrow 
-\frac{\mathit{\ssa} \eta_{1}}{\mathit{\ell}_{1}}-\frac{\mathit{\ssb} \eta_{2}}{\mathit{\ell}_{2}}+p_{1}, \quad
\frac{d^{2}}{d t^{2}}x_{1}\! \left(t \right) \leftrightarrow 
-\frac{\left(\eta_{1}+\eta_{2}\right) \left(\mathit{\cca} \eta_{1}+\mathit{\ccb} \eta_{2}+\mathit{\ell}_{1} p_{2}\right)}{\mathit{\ell}_{1}^{3}}, \\
\frac{d}{d t}x_{2}\! \left(t \right) & \leftrightarrow
\frac{\eta_{1} \mathit{\cca}}{\mathit{\ell}_{1}}+\frac{\eta_{2} \mathit{\ccb}}{\mathit{\ell}_{1}}+p_{2}, \quad \frac{d^{2}}{d t^{2}}x_{2}\! \left(t \right) \leftrightarrow 
\frac{\left(\eta_{1}+\eta_{2}\right) \left(-\mathit{\ssa} \eta_{1}-\mathit{\ssb} \eta_{2}+\mathit{\ell}_{1} p_{1}\right)}{\mathit{\ell}_{1}^{3}}.
\end{align*}
As a result, 
$$
\left(\frac{d}{d t}x_{1}\! \left(t \right)\right) \left(\frac{d^{2}}{d t^{2}}x_{2}\! \left(t \right)\right)-\left(\frac{d^{2}}{d t^{2}}x_{1}\! \left(t \right)\right) \left(\frac{d}{d t}x_{2}\! \left(t \right)\right)
\leftrightarrow   \kappa poly_x =
\frac{\eta_{1}}{\ell_{1}^{2}}+\frac{\eta_{2}}{\ell_{1}^{2}} \mod I.
$$
We use $\kappa_x$ to denote the curvature for the curve $x =(x_1(t), x_2(t))$, $\kappa poly_x$  the associated polynomial in terms of the variables {\it vars}. Note that this agrees with \ref{eq:curv1} .

As a brief aside, one can compute {\it mutatis mutandis} the curvatures of the curves $y_1$ and $y_2$  using this same technique, whereas computing by hand is somewhat messy. Here is the polynomial associated with the $y_1$ curvature, the $y_2$ formula is similar:

$$
\kappa poly_ {y _1 }=
\frac{2 \mathit{\cca} \mathit{\ccb} \eta_{1}}{\mathit{\ell}_{1}^{2}}+\frac{2 \mathit{\cca} \mathit{p_2}}{\mathit{\ell}_{1}}+\frac{2 \eta_{1} \mathit{\ssa} \mathit{\ssb}}{\mathit{\ell}_{1}^{2}}-\frac{2 \mathit{p_1} \mathit{\ssa}}{\mathit{\ell}_{1}}-\frac{\eta_{1}}{\mathit{\ell}_{1}^{2}}+\frac{\eta_{2}}{\mathit{\ell}_{1}^{2}}.
$$
For computational purposes, this expression is the one we wish to work with.  The reader is invited to manipulate the associated functions, and come up with a simpler expression.

We return to the elastica calculation, the case when $\ell_1=\ell_2$.  
We wish to show that there exists a constant $A$ such that $\kappa = \kappa_x$
satisfies 
\begin{equation} \label{eq:elasticade}
\frac{d^{2}}{d t^{2}}\kappa \! \left(t \right) = 
-\frac{\kappa \! \left(t \right)^{3}}{2}-A \kappa \! \left(t \right).
\end{equation}
Define  $\kappa_0, \kappa_1, \kappa_2$ as the polynomials associated with $\kappa, \kappa', \kappa''$ respectively.  At each stage,
differentiate a function, substitute the differential equations system (\ref{eq:bigdes}), convert to a polynomial, and then compute and reduce via the Gr\"obner basis {\it GPB}.  Substitute the derivatives in \ref{eq:elasticade}; one obtains
an equation for $A$, with solution which has the corresponding polynomial:

$$
Apoly = -\frac{\mathit{\cca} p_{2} \eta_{1}}{\mathit{\ell}_{1}^{3}}-\frac{\mathit{\ccb} \eta_{2} p_{2}}{\mathit{\ell}_{1}^{3}}+\frac{\eta_{2} p_{1} \mathit{\ssb}}{\mathit{\ell}_{1}^{3}}-\frac{p_{1}^{2}}{\mathit{\ell}_{1}^{2}}-\frac{p_{2}^{2}}{\mathit{\ell}_{1}^{2}}+\frac{p_{1} \mathit{\ssa} \eta_{1}}{\mathit{\ell}_{1}^{3}}-\frac{\eta_{2}^{2}}{2 \mathit{\ell}_{1}^{4}}-\frac{\eta_{2} \eta_{1}}{\mathit{\ell}_{1}^{4}}-\frac{\eta_{1}^{2}}{2 \mathit{\ell}_{1}^{4}}.
$$
To show that $A$ is actually constant, differentiate and substitute the equations (\ref{eq:bigdes}); the expression reduces to $0$, as desired.

Essentially all results in this paper can be proved directly by hand in the equal length case, and proofs are given.  Our results  for unequal lengths depend heavily on the Gr\"obner basis reduction tools.

\section{B\"acklund transformations of 1-soliton curves to 2-soliton curves} \label{sect:Ba}

In this section, we recall an elementary contruction in curve theory.  For our purposes, we only need the contruction for planar curves.  We seek to answer the following question:  given a curve $\gamma$, how can we construct a new curve 
$\tilde \gamma$, with a point-to-point correspondence between the new curves, such that 
\begin{enumerate}
\item the segments representing the point-to-point correspondence have {\it constant} length;
\item the arc length parameterization is preserved via the correspondence?  
\end{enumerate}
The construction has various names:  ``bicycle correspondence", B\"acklund transformation (see \cite{BLPT,Ta}, and \cite{RS}).  

Specifically, given $\gamma(t)$ with arc length parameter $t$,  we construct a new curve
$\tilde \gamma(t)$  via
 $$
 \tilde \gamma (t) =\gamma(t)  + L\, ( \cos \beta(t)\, T(t) + \sin\beta(t)\, N(t)), 
 $$
  where $L$ is constant and $(T,N)$ is the Frenet frame along $\gamma$. 
  
 Note that condition 1) is automatically satisfied, with segment length $L$, but that $t$ is initially only known to be the arc length parameter along $\gamma$, not $\tilde \gamma$.  We ask when the $t$ parameter in $\tilde \gamma$ is actually arc length along the second curve. Of course, this means that
$ {\tilde \gamma}' \cdot  {\tilde \gamma}' -1 = 0,$ and implies the following restriction on $\beta$:
$$
\left(\mathit{\beta'(t)} +\kappa \! \left(t \right)\right) \left(L \mathit{\beta'} +L \kappa \! \left(t \right)-2 \sin \! \left(\beta \! \left(t \right)\right)\right)=0,
$$
where $\kappa(t)$ is the curvature of $\gamma$. 

If the first factor vanishes, then $\beta' = -\kappa$,  which corresponds to the displacement field 
$L(\cos \beta T + \sin \beta N)$ being constant along $\gamma$.  Otherwise, we have the first order differential equation for $\beta$:
$$
\frac{d}{d t}\beta \! \left(t \right) = 
-\kappa \! \left(t \right)+\frac{2 \sin \! \left(\beta \! \left(t \right)\right)}{L}
$$
(this is the ``bicycle" ODE, compare with Remark \ref{rmk:bic}).
Using the differential equation, one can easily compute that:
\begin{equation*}
\begin{aligned}
\tilde T &= 
T \! \left(t \right) \cos \! \left(2 \beta \! \left(t \right)\right)+
N \! \left(t \right) \sin \! \left(2 \beta \! \left(t \right)\right); \\
\tilde N &=
-N \! \left(t \right) \cos \! \left(2 \beta \! \left(t \right)\right)+
T \! \left(t \right) \sin \! \left(2 \beta \! \left(t \right)\right); \\
\tilde \kappa &= 
\kappa \! \left(t \right)-\frac{4 \sin \! \left(\beta \! \left(t \right)\right)}{L}.
\end{aligned}
\end{equation*}

If $\gamma$ is an elastic curve, we can the describe the geometry of its B\"acklund transformation $\tilde \gamma$.

\begin{thm} \label{thm:one2two}
If $\gamma$ is an elastic curve, then $\tilde \gamma$ is a 2-soliton.
\end{thm}

\bp
We are given the  equations 
\begin{equation}\label{eq:bdes}
\begin{aligned}
\left(\frac{d}{d t}\kappa \! \left(t \right)\right)^{2} &= 
-\kappa \! \left(t \right)^{2} A -\frac{\kappa \! \left(t \right)^{4}}{4}-B, \quad
\frac{d^{2}}{d t^{2}}\kappa \! \left(t \right) = 
-\kappa \! \left(t \right) A -\frac{\kappa \! \left(t \right)^{3}}{2}, \\
\tilde \kappa &= 
\kappa \! \left(x \right)+\frac{4 \sin \! \left(\beta \! \left(x \right)\right)}{L}, \quad
\frac{d}{d t}\beta \! \left(t \right)  = 
-\kappa \! \left(t \right)-\frac{2 \sin \! \left(\beta \! \left(t \right)\right)}{L},
\end{aligned}
\end{equation}
and we wish to conclude that $\tilde \kappa$ satisfies, for appropriate constants $c_1, c_2$, the 
equation
\begin{equation}\label{eq:eul4}
\begin{aligned}
\kappa \! \left(t \right) c_{1}+\left(-\frac{d^{2}}{d t^{2}}\kappa \! \left(t \right)-\frac{\kappa \! \left(t \right)^{3}}{2}\right) c_{2}&+\frac{5 \kappa \! \left(t \right)^{2} \left(\frac{d^{2}}{d t^{2}}\kappa \! \left(t \right)\right)}{2}+\frac{5 \kappa \! \left(t \right) \left(\frac{d}{d t}\kappa \! \left(t \right)\right)^{2}}{2}\\
&+\frac{3 \kappa \! \left(t \right)^{5}}{8}+\frac{d^{4}}{d t^{4}}\kappa \! \left(t \right)=0,
\end{aligned}
\end{equation}
of course, with $\kappa$ replaced by $\tilde \kappa$.

We define $\kappa_0 = \tilde \kappa, \kappa_1 = \kappa_0' , \dots \kappa_4=\kappa_3'$.  At each stage after the differentiation, we simplify the expression using (\ref{eq:bdes}); the resulting expressions  are polynomial in $\cos\beta, \sin\beta$ with no $\beta$
derivatives, and only involve polynomial expressions in $\kappa$ and, at most, a term {\it linear} in $\kappa'$.

The coefficient of $\kappa'$ in the above expression is
$$
\frac{4 \left(\kappa  L +4 \sin \! \left(\beta \right)\right) \left(\left(\frac{B}{8}
+\frac{c_{1}}{4}\right) L^{2}+A \right)}{L^{3}},
$$
hence  
$$c_{1} = -\frac{B}{2}-\frac{4 A}{L^{2}}.$$
 Substituting this value back in our equation, we obtain
$$c_{2} = -A +\frac{4}{L^{2}},$$
as needed.
\ep


We saw in Section \ref{sect:ela} that, in the equal length case, the path of the ``front wheel" $x$ is an elastica.
We can now describe the geometry of the trajectories of  the ``rear points" $y_1, y_2$.

By construction, the paths of the rear points $y_1, y_2$ are in B\"acklund correspondence with the path of $x$ (an elastica)
with  lengths equal to $2 \ell$.  Hence we can conclude:

\begin{cor}\label{cor:equal2sol}
The paths of $y_1, y_2$ are 2-soliton curves. 
\end{cor}

\bp
By Theorem \ref{thm:one2two}, the curvatures $\kappa_{y_1}, \kappa_{y_2}$, satisfy (\ref{eq:eul4})
with coefficients $c_1, c_2$:
$$c_{1} = -\frac{B}{2}-\frac{4 A}{L^{2}} , \quad c_{2} = -A +\frac{4}{L^{2}},$$
where $L = 2 \ell$, and $A,B$ are related as in Lemma \ref{lm:AB}.
\ep


\section{Examples}\label{sect:examples}
\begin{example}\label{ex:circle}
{\rm
Consider the special case when $x$ traverses a unit circle:  $\gamma(t)=t$. In this case, equations (\ref{eq:abg}) can be solved exactly.

One type of solutions is when $\alpha_2=t-\pi/2$ and $\alpha_1(t)$ is any function satisfying  $\dot \alpha_1 = {\sin(t-\alpha_1)}$. This solution corresponds, via lifting, to the special solution of the 1-linkage problem when the rear wheel of the bicycle is fixed at the center of the unit circle. Of course, by symmetry, there is a solution with $\alpha_1=t-\pi/2$ and $\alpha_2(t)$ is any function satisfying  $\dot \alpha_2 = {\sin(t-\alpha_2)}$.

  A calculation that we omit shows that there are no other solutions to (\ref{eq:abg}).
}
\end{example} 

\begin{figure}[ht] 
\centering
\includegraphics[width=.48\textwidth]{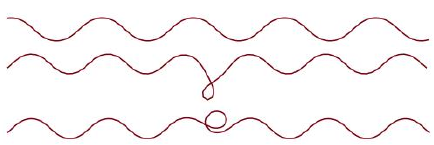}
\includegraphics[width=.48\textwidth]{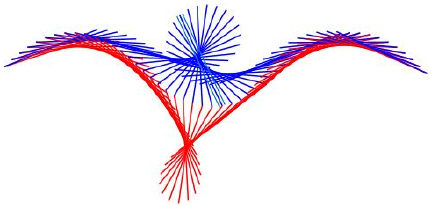}
\includegraphics[width=.48\textwidth]{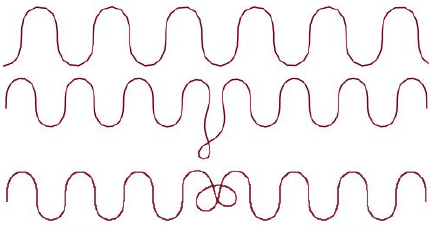}\qquad
\includegraphics[height=1.5in]{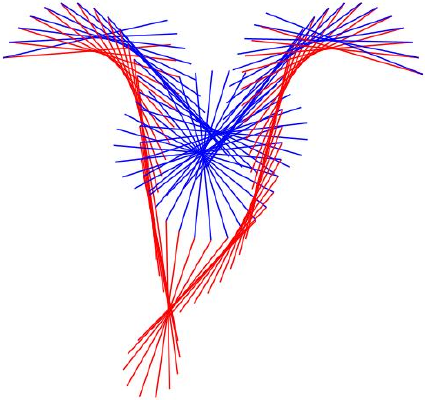}
\caption{On the left, the trajectories of points $x, y_1$, and $y_2$ are shown, and on the right, the respective motion of the 2-linkage whose two links are colored red and blue. The values of $k$ are $0.1$ and $0.707$, respectively.}	
\label{Ex1}
\end{figure}

\begin{example}\label{ex:inflection}
{\rm
Here,  we generate initial conditions for equation (\ref{eq:bigdes}) which produce various inflectional elasticae. We assume that $\ell_1 = \ell_2 = 1$. 
For convenience, the variables used here will (by slight abuse of notation)
represent functions values at $0$.  For example, $\eta_1$ represents $\eta_1(0))$, $\cca$ represents $\cos(\alpha_1(0)$, etc.

It is well known that the function $\nu= 2 k  \, \text{cn}(t,k), \, 0 < k < 1$, is the formula for the curvature of  inflectional elastica 
(the limiting cases of $k = 0,\, 1$  correspond to the line and Euler soliton respectively). Thus, we have $\nu(0) = 2k$ and $\nu'(0) = 0$.
Also, we have 
$$\nu'' + \frac{\nu^3}{2} + A_{\nu}\nu = 0,
$$
 where $A_{\nu} = -2 k^{2}+1$.

\begin{figure}[ht] 
\centering
\includegraphics[height=2.5in]{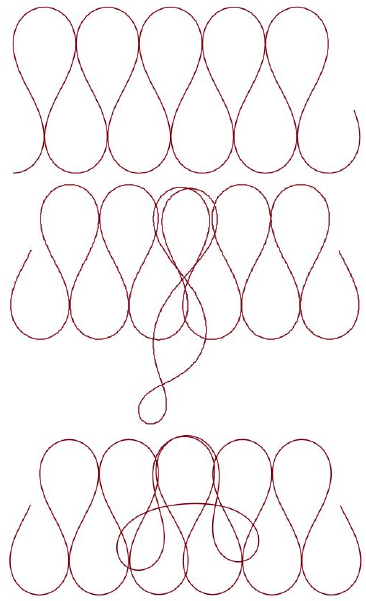}\qquad
\includegraphics[height=2.3in]{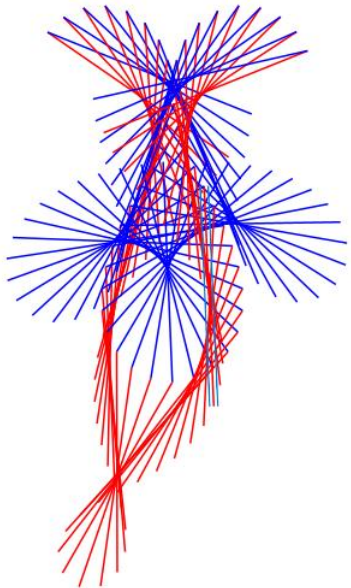}
\includegraphics[height=1.8in]{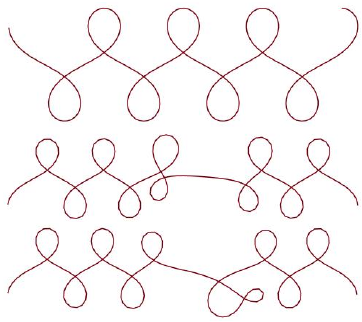}\qquad
\includegraphics[height=1.8in]{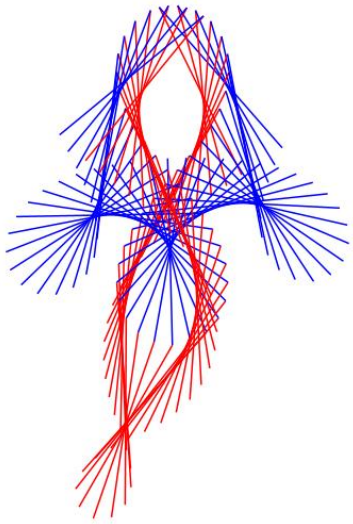}
\caption{On the left, the trajectories of points $x, y_1$, and $y_2$ are shown, and on the right, the respective motion of the 2-linkage.
The values of $k$ are $0.854$ and $0.95$, respectively.}	
\label{Ex2}
\end{figure}

On the other hand, in Section \ref{sect:comp} we saw that the values of  $\kappa_x$ and its derivative at $0$ are given by 
$\eta_{1}+\eta_{2}$  and
$$
   {\cca} \eta_{1} p_{1}+   {\ccb} \eta_{2} p_{1}+   {\ssa} \eta_{1} p_{2}+   {\ssb} \eta_{2} p_{2},
$$
respectively.  In addition, we know that $\kappa_x$ satisfies the differential equation
$$
\kappa'' + \frac{\kappa^3}{2}+ A_L \kappa= 0,
$$
 where $A_L$ is given by ($L$ stands for linkage)
\begin{align*}
A_L =&-\cca p_2 \eta_1 - \ccb \eta_2 p_2 + \eta_2 p_1 \ssb + p_1 \ssa \eta_1  \\
&- \frac{1}{2} \eta_2^2 - \eta_2 \eta_1 - p_1^2 - p_2^2 - \frac{1}{2} \eta_1^2. 
\end{align*}
If we can find initial conditions satisfying
\begin{align*}
&\eta_{1}+\eta_{2} = 2 k,  \\
& {\cca} \eta_{1} p_{1}+   {\ccb} \eta_{2} p_{1}+   {\ssa} \eta_{1} p_{2}+   {\ssb} \eta_{2} p_{2}=0 , \\
&A_L = A_\nu
\end{align*}
with the additional condition that $H- \frac{1}{2} = 0$ (so that we have unit speed parameterization of the $(x_1, x_2)$ curve),
then the functions $\nu$ and $\kappa_x$ satisfy the same differential equation with same initial conditions; and, by the standard theorem
of existence and uniqueness of solutions to ordinary differential equations (e.g.,\cite{CL}), the two functions are equal.

\begin{figure}[ht] 
\centering
\includegraphics[height=1.3in]{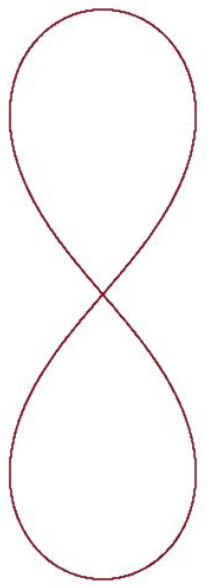}\qquad
\includegraphics[height=1.5in]{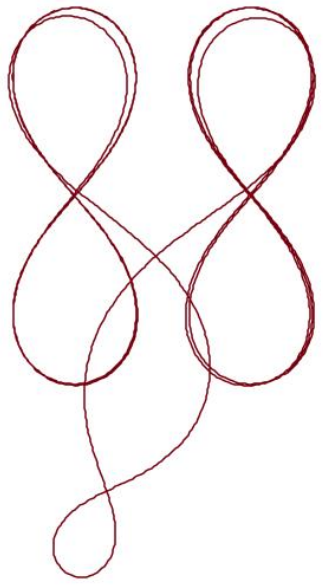}\quad
\includegraphics[height=1.3in]{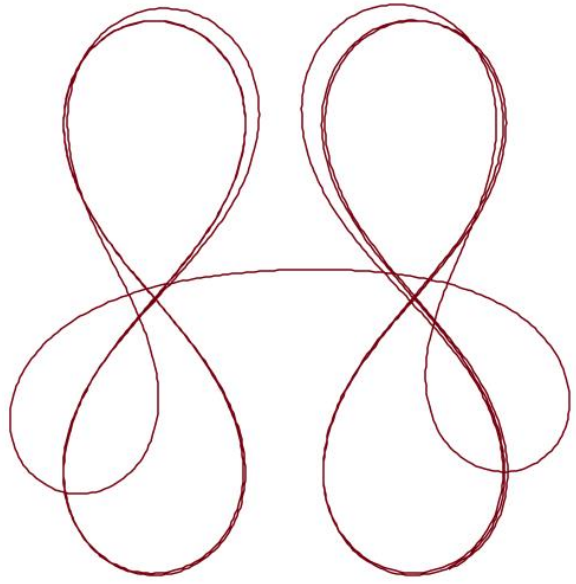}\quad
\includegraphics[height=1.7in]{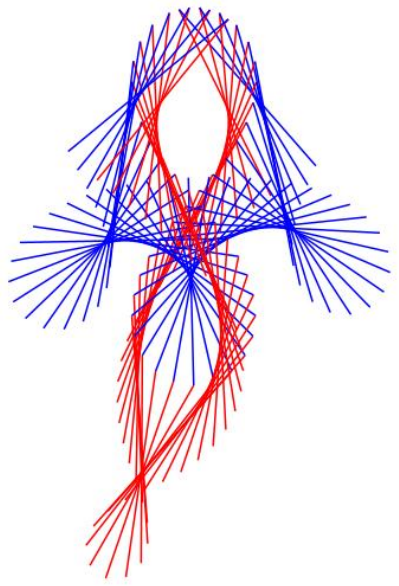}
\caption{From left to right: the trajectories of points $x, y_1,y_2$ and the respective motion of the 2-linkage. Here $k=0.909$.}	
\label{Ex3}
\end{figure}

There are many solutions to this set of equations; by judicious choice of $\alpha_1, p_1, p_2$, we obtain the solution:
\begin{align*}
& \alpha_{1}\,     (0     ) = 
\arctan \,     \left(\frac{4}{3}\right) , \, \alpha_{2}\,     (0     ) = 
\arctan \,     \left(\frac{4 k^{2}-10 k +4}{3 k^{2}-3}     \right), 
\eta_{1}\,     (0     ) = \frac{5 k^{2}-5}{5 k -4}, \\
&\eta_{2}\,     (0     ) = \frac{5 k^{2}-8 k +5}{5 k -4}, 
p_{1}\,     (0     ) = 1, p_{2}\,     (0     ) = 0, 
x_{1}\,     (0     ) = 0, x_{2}\,     (0     ) = 0.
\end{align*}

Inserting values $k = .1, \, .707, \, .854, \, .95$, and $.909$,  one obtains the initial conditions which produce  various different types of 
inflectional elasticae, see Figure \ref{Ex1},\ \ref{Ex2}, and \ref{Ex3}.

Note that, in the last case (Figure \ref{Ex3}), the trajectory of point $x$ is a closed curve, but the trajectories of points $y_1$ and $y_2$ are not. We do not know whether there exists a closed tricycling geodesic whose $x$-projection is an eight-shaped elastica.

}
\end{example}

\section{The case of unequal length line segments and the  re-appearance of 2-soliton curves} \label{sect:un2sol}

\begin{figure}[ht] 
\centering
\includegraphics[width=.6\textwidth]{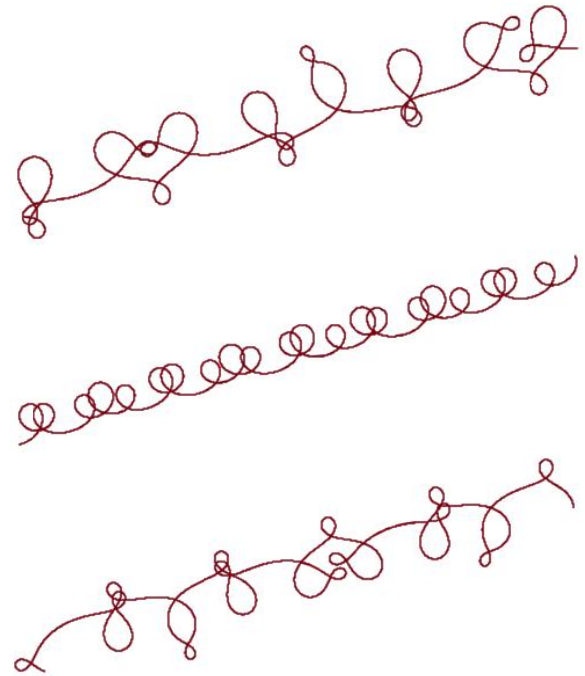}
\caption{From top to bottom: the trajectories of points $x, y_1$, and $y_2$.}	
\label{diff}
\end{figure}

Before making any hypotheses about the geometry of the paths of $x, y_1, y_2$, we looked at some plots of these paths -- see Figure \ref{diff}.

 It is immediately obvious that {\it none} of the three paths are elastica.  However, given the results of the last section, it is plausible that these curves are 2-soliton curves.
In this section, we prove that this is the case.

\begin{thm} \label{thm:noneq}
If $\ell_1 \neq \ell_2$, all three planar projections of a tricycling geodesic are 2-solitons.
\end{thm}

\bp
We begin with the $x$-curve.  We know that its curvature is given by
$$\kappa_x = \frac{\eta_{2}\! \left(t \right)}{\mathit{\ell_2}^{2}}+\frac{\eta_{1}\! \left(t \right)}
{\mathit{\ell_1}^{2}},$$
and we wish to show, as in the equal length case, that $\kappa_x$ satisfies the differential equation 
\begin{equation}\label{eq:eul4b}
\begin{aligned}
a \kappa \! \left(t \right)+b \left(-\frac{d^{2}}{d t^{2}}\kappa \! \left(t \right)-\frac{\kappa \! \left(t \right)^{3}}{2}\right)+\frac{5 \kappa \! \left(t \right)^{2} \left(\frac{d^{2}}{d t^{2}}\kappa \! \left(t \right)\right)}{2}+\frac{5 \kappa \! \left(t \right) \left(\frac{d}{d t}\kappa \! \left(t \right)\right)^{2}}{2}\\
+\frac{3 \kappa \! \left(t \right)^{5}}{8}+\frac{d^{4}}{d t^{4}}\kappa \! \left(t \right) =0
\end{aligned}
\end{equation}
for appropriate $a, b$.  

As before, we can set $\kappa_0 = \kappa_x, \kappa_1 = \kappa_0', \dots \kappa_4 = \kappa_3'$, again simplifying  the expressions at each step using (\ref{eq:bigdes}) and algebraic reduction,  then replacing the derivatives of $\kappa$ by the 
appropriate $\kappa_i$.
However, this will simply give us one condition on the two unknowns $a,b$.  

We remedy this situation by differentiating
\ref{eq:eul4b} to obtain a {\it fifth} order equation for $\kappa_x$, requiring an additional computation of $\kappa_5 = \kappa_4'$.
After algebraic reduction, one obtains the polynomial representation for $b$:
\begin{align*}
bpoly&=\frac{\mathit{\cca} p_{2} \eta_{1}}{\ell_{1} \ell_{2}^{2}}+\frac{\mathit{\ccb} \eta_{2} p_{2}}{\ell_{1}^{2} \ell_{2}}
-\frac{\eta_{2} p_{1} \mathit{\ssb}}{\ell_{1}^{2} \ell_{2}}+\frac{p_{1}^{2}}{2 \ell_{2}^{2}}+\frac{p_{2}^{2}}{2 \ell_{2}^{2}}\\&-\frac{p_{1} \mathit{\ssa} \eta_{1}}{\ell_{1} \ell_{2}^{2}}+\frac{p_{1}^{2}}{2 \ell_{1}^{2}}+\frac{p_{2}^{2}}{2 \ell_{1}^{2}}+\frac{\eta_{2}^{2}}{2 \ell_{1}^{2} \ell_{2}^{2}}+\frac{\eta_{1} \eta_{2}}{\ell_{1}^{2} \ell_{2}^{2}}+\frac{1}{2 \ell_{2}^{2}}+\frac{1}{2 \ell_{1}^{2}}+\frac{\eta_{1}^{2}}{2 \ell_{1}^{2} \ell_{2}^{2}}
\end{align*}
and a considerably longer expression for $a$. One can check directly that $a,b$ are constants of motion.

In anticipation of our discussion of the complete integrability of (\ref{eq:bigdes}), we note that $H, p_1, p_2$ are constants
of motion.  Our configuration space is four-dimensional, so we could anticipate an addtional constant of motion, but here we discover two of them, $b$ and $a$.  We hypothesize that $a$ is functionally dependent on $b$.  

Thus, we augment
our polynomial basis $PB_2 = PB \, \bigcup \, \{bpoly - \lambda\}$, corresponding to a larger ideal
$I_2$, and use the appropriate tool to compute a
new Gr\"obner basis $GPB_2$ for $I_2$.  One obtains
$$
apoly = 
\frac{p_{1}^{2}}{2 \ell_{1}^{2} \ell_{2}^{2}}+\frac{p_{2}^{2}}{2 \ell_{1}^{2} \ell_{2}^{2}}-\frac{\lambda^{2}}{2}+\frac{\lambda}{\ell_{2}^{2}}+\frac{\lambda}{\ell_{1}^{2}}-\frac{1}{2 \ell_{2}^{4}}-\frac{1}{2 \ell_{1}^{2} \ell_{2}^{2}}-\frac{1}{2 \ell_{1}^{4}} \mod I_2.
$$
Hence, by replacing $\lambda$ with $b$ in the previous formula, we see that $a$ is in fact a quadratic function of $b$. 

One can also check, using the same sort of computations, that $\kappa_{y_1}, \kappa_{y_2}$ both satisfy \ref{eq:eul4b} with the {\it same} values of $a$ and $b$ just computed.  
\ep
In the discussion of complete integrability which follows, it will be convenient to use an alternative notation $\tilde G = b$.

\section{Complete integrability of the 2-linkage geodesic equations} \label{sect:int}
In this section we show that the geodesic equations of the 2-linkage (\ref{eq:Ham}) are completely integrable.
We begin with the case of equal lengths.

Consider the cotangent bundle $T^*{\mathcal C}$ with its standard symplectic form 
$$
dx_1\wedge dp_1 + dx_2\wedge dp_2 + d\alpha_1\wedge d\eta_1 + d\alpha_2\wedge d\eta_2.
$$
We have four integrals of motion
\begin{equation} \label{eq:integr}
\begin{aligned}
&H =  \frac{1}{2} [(p_1-\eta_1 \sin\alpha_1 -\eta_2 \sin\alpha_2)^2 + (p_2+\eta_1 \cos\alpha_1 +\eta_2 \cos\alpha_2)^2], p_1, p_2, \\
&G=\frac{\ddot \kappa}{\kappa} + \frac{1}{2} \kappa^2 = p_1 (p_1-\eta_1 \sin\alpha_1 -\eta_2 \sin\alpha_2) + p_2 (p_2+\eta_1 \cos\alpha_1 +\eta_2 \cos\alpha_2)\\
 &+ \frac{1}{2}(\eta_1+\eta_2)^2,
\end{aligned}
\end{equation}
where we used (\ref{eq:curv1}) and (\ref{eq:doubleprime}).

\begin{prop} \label{prop:comm}
These integrals Poisson commute and they are almost everywhere functionally independent.
\end{prop}

\bp
It is obvious that $p_1$ and $p_2$ are in involution with everything. That $\{H,G\}=0$
follows from the fact that $G$ is constant along the trajectories of the Hamiltonian field with the Hamiltonian function $H$.

For the second statement, let 
$$
u_1=\eta_1\cos\alpha_1+\eta_2\cos\alpha_2,\ u_2=\eta_1\sin\alpha_1+\eta_2\sin\alpha_2.
$$
The change of variables
$$
(\alpha_1,\alpha_2,\eta_1,\eta_2) \mapsto (u_1,u_2,\eta_1,\eta_2)
$$
has the Jacobian $\eta_1\eta_2\sin(\alpha_2-\alpha_1)$. Assume that this Jacobian doesn't vanish.
Then 
$$
H=\frac{(p_1-v)^2+(p_2+u)^2}{2},\ G=p_1^2+p_2^2 + p_2 u_1 - p_1 u_2 +\frac{(\eta_1+\eta_2)^2}{2}.
$$
Since $p_1$ and $p_2$ are integrals, we need to see when $\nabla H$ and $\nabla G$ (as functions of $u_1,u_2,\eta_1,\eta_2$) are linearly dependent. 
The relevant $2\times 4$ matrix is
$$
	\begin{pmatrix} 
	p_2+u_1&-p_1+u_2&0&0\\
	p_2&-p_1&\eta_1+\eta_2&\eta_1+\eta_2\\
	\end{pmatrix}.
$$

If $\kappa=\eta_1+\eta_2 \neq 0$, then this matrix does not have full rank if and only if $u_2=p_1, u_1=-p_2$, but this is impossible because then one would have $H=0$. 

If $\kappa=0$, then the matrix does not have full rank if and only if  $p_1 u_1 + p_2 u_2=0$, that is, 
$$
p_1 (\eta_1\cos\alpha_1+\eta_2\cos\alpha_2) + p_2 (\eta_1\sin\alpha_1+\eta_2\sin\alpha_2)=0.
$$
This is a hypersurface in the phase space.

To see the meaning of this condition, substitute $\eta_1$ and $\eta_2$ from (\ref{eq:elim}) to arrive (as always, applying trigonometric identities) at $p_2 \cos\gamma - p_1 \sin\gamma=0$ or, in view of (\ref{eq:prime}), $\dot \kappa =0$. 
\ep

We briefly discuss the case of unequal lengths.  Here, we introduce a slightly different change of variables (which fails to work in the equal length case).

Define
\begin{align*}
\lambda_{1}& = 
\frac{\cos \! \left(\alpha_{1}\right) \eta_{1}}{\mathit{\ell}_{1}}+\frac{\cos \! \left(\alpha_{2}\right) \eta_{2}}{\mathit{\ell}_{2}}
, \
\lambda_{2} = 
\frac{\sin \! \left(\alpha_{1}\right) \eta_{1}}{\mathit{\ell}_{1}}+\frac{\sin \! \left(\alpha_{2}\right) \eta_{2}}{\mathit{\ell}_{2}}
, \\
\mu_{1} &= 
\frac{\eta_{2} \sin \! \left(\alpha_{2}\right)}{\mathit{\ell}_{1}}+\frac{\sin \! \left(\alpha_{1}\right) \eta_{1}}{\mathit{\ell}_{2}}
,\ \mu_{2} = \eta_{2}+\eta_{1}.
\end{align*}
One can check that the determinant of the Jacobian of the transformation 
$(\alpha_{1}, \alpha_{2}, \eta_{1}, \eta_{2}) \rightarrow (\lambda_{1}, \lambda_{2}, \mu_{1}, \mu_{2})$
is
$$
\frac{\eta_{1} \eta_{2} \left(\cos \! \left(\alpha_{1}\right) \mathit{\ell}_{1}-\cos \! \left(\alpha_{2}\right) \mathit{\ell}_{2}\right) \left(\mathit{\ell}_{1}^{2}-\mathit{\ell}_{2}^{2}\right)}{\mathit{\ell}_{1}^{3} \mathit{\ell}_{2}^{3}} .
$$

We can express our two Hamiltonians $H$ and $\tilde G$ in terms of these new variables (WLOG, we  henceforth make the assumption that $p_2 = 0$):
$$
H = 
\frac{1}{2} \lambda_{1}^{2}-\lambda_{1} p_{1}+\frac{1}{2} \lambda_{2}^{2}+\frac{1}{2} p_{1}^{2},
\quad
\tilde G = 
\frac{p_{1}^{2}}{2 \mathit{\ell}_{2}^{2}}-\frac{\mu_{1} p_{1}}{\mathit{\ell}_{2} \mathit{\ell}_{1}}+\frac{p_{1}^{2}}{2 \mathit{\ell}_{1}^{2}}+\frac{1}{2 \mathit{\ell}_{2}^{2}}+\frac{1}{2 \mathit{\ell}_{1}^{2}}+\frac{\mu_{2}^{2}}{2 \mathit{\ell}_{2}^{2} \mathit{\ell}_{1}^{2}}.
$$
To prove independence of the two functions $H$ and $\tilde G$, we need to show that the matrix $[\nabla H, \nabla G]$ has full rank; but this
matrix is easily computed to be:
$$
\left[\begin{array}{cccc}
\lambda_1 -\mathit{p_1}  & \lambda_2  & 0 & 0 
\\
 0 & 0 & -\frac{\mathit{p_1}}{\mathit{\ell_2} \mathit{\ell_1}} & \frac{\mu_2}{\mathit{\ell_2}^{2} \mathit{\ell_1}^{2}} 
\end{array}\right].
$$
Ignoring the trivial case $p_1 = 0$ (where the $x$-path is a circle), we see that we have dependence only if both $\lambda_2$ and
$\lambda_1 - p_1$ are simultaneously equal to zero.

\section{Geometric interpretation of the constants of motion $ G, \tilde G$.  The planar filament equation} \label{sect:fil}
Recall that for the equal length case, we have a constant of motion $ G$, whereas we had a constant of motion $\tilde G$ 
in the unequal case.  One can check that if  $\ell_1=\ell_2  = \ell$, then $\tilde G - G= \frac{1}{\ell^2}$.  In this section we discuss
a geometric interpretation of $ G, \tilde G$ in their respective cases.

First, assume that $\ell_1=\ell_2=1$.

Fix a point of $Q\in T^*{\mathcal C}$, and consider the flow lines of the commuting Hamiltonian vector fields
 $X_H$ and $X_{G}$ through point $Q$, where $H$ and $G$ are as in (\ref{eq:integr}). 
 
 Project these curves to the plane. 
The first curve is the trajectory $x(t)$, an elastica. Let $(T(t),N(t))$ be the Frenet frame along this curve. 
The tangent vectors to the second curve provide a vector field $\xi$ along this elastica. This vector field can be 
written as $\xi=uT+vN$, where $u$ and $v$ are functions of $t$. 

\begin{thm} \label{thm:fil}
One has
$$
u(t)=1+g-\frac{1}{2}\kappa^2(t),\ v(t)=-\dot \kappa(t),
$$
where $\kappa(t)$ is the curvature of the curve $x(t)$, and $g$ is the value of the integral $G$ at point $Q$.
\end{thm}

\bp
The projection of the vector  $X_{\tilde G}(Q)$ on the plane is the vector 
\begin{equation*}
\begin{aligned}
\xi=(G_{p_1},G_{p_2})&=(2p_1-\eta_1\sin\alpha_1-\eta_2\sin\alpha_2, 2p_2+\eta_1\cos\alpha_1+\eta_2\cos\alpha_2)\\
&=(p_1+\cos\gamma,p_2+\sin\gamma),
\end{aligned}
\end{equation*}
where the last equality is due to (\ref{eq:fnct}). One has
$$
T=(\cos\gamma,\sin\gamma), N=(\sin\gamma,\cos\gamma),
$$
and the linear equation $\xi=uT+vN$ has the solution
$$
u=1+p_1\cos\gamma+p_2\sin\gamma, v=-p_1\sin\gamma+p_2\cos\gamma.
$$
Use equations (\ref{eq:doubleprime}) and (\ref{eq:prime}) and the definition of $G$ to conclude:
$$
u = 1 + \frac{\ddot \kappa}{\kappa}= 1+g-\frac{1}{2}\kappa^2(t),\ v(t)=-\dot \kappa(t),
$$ 
as claimed.
\ep

Thus the vector field $\xi$ is a linear combination of the unit tangent field $T$ and the planar filament vector field $\frac{1}{2}\kappa^2 T + \dot \kappa N$, see \cite{LP}. Note that the flow of this field takes the elastic curves to elastic curves, in agreement with the fact that the flows of the Hamiltonians $H$ and $G$ commute.

If $\ell_1 \neq \ell_2$,  we have the analogous statement, proved by the computer algebra techniques already discussed.

\begin{thm} \label{thm:filb}
One has
$$
u(t) =\tilde g -\frac{1}{2} \kappa^2(t),\ v(t)=-\dot \kappa(t),
$$
where $\kappa(t)$ is the curvature of the curve $x(t)$, and $\tilde g$ is the value of the integral $\tilde G$ at point $Q$.
\end{thm}


\section{Concluding remarks} \label{sect:concl}

Inspired by papers \cite{2D, BJT}, we have obtained a 2-linkage model which exhibits the fundamental observation
in those references: the vertices of the moving linkage sweep out curves of interesting geometry.
Amusingly, our problem, involving a completely non-integrable distribution,  is completely integrable.  

\begin{figure}[ht] 
\centering
\includegraphics[width=1\textwidth]{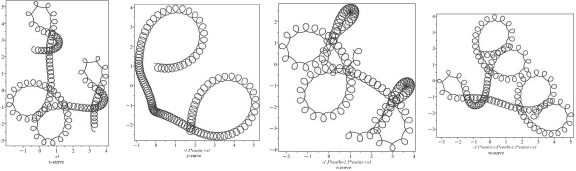}
\caption{The trajectories of four vertices of a ``linear" 3-linkage, the case of three unequal lengths.}	
\label{lin3diff}
\end{figure}

\begin{figure}[ht] 
\centering
\includegraphics[width=.8\textwidth]{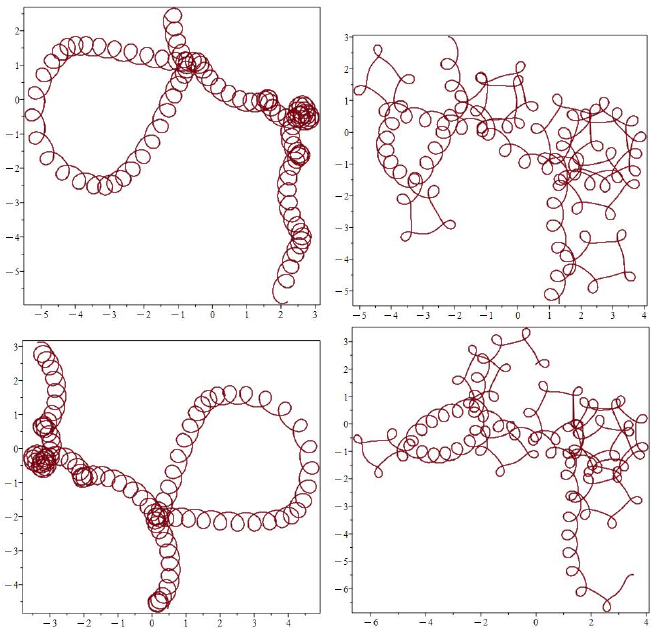}
\caption{The trajectories of four vertices of a ``spokewheel" 3-linkage, the case of three unequal lengths. The first trajectory is of that of the middle point.}	
\label{wheel}
\end{figure}

It is clear that there is much more worth studying related to this problem.  Here are some open questions:
\begin{itemize}
\item Elastic curves have been studied in the context of spherical and hyperbolic geometry, see \cite{LS2}; additionally,
analogues of the planar filament flow exist in those geometries as well, see \cite{LP2}.  This suggests that our 2-linkage 
problem could be amenable to study is these contexts too.
\item In this paper, we have by no means solved the geodesic equations, we have only shown a certain striking
property of the projection of the solutions to the plane.  Can these equations be solved exactly?
\item Paper \cite{BJT}  extends the work in \cite{2D} to 1-linkages in $\R^n$. We expect  the results of our paper to be extendable to higher dimensional cases as well.
\item Isospectrality is a nearly universal property of integrable systems, although this can take some work to uncover.
To give one example, the problem of integrability of geodesics on the ellipsoid goes back to Jacobi, but the isospectral
framework for the problem was discovered by Moser \cite{Mos}.  Does such a framework appear in the case of
our problem?
\item The starting point of our work was a different sub-Riemmanian problem, concerning the geodesic motion of a 
2-linkage $xyz$ subject to the following non-holonomic constraint: the velocity of point $y$ is aligned with the segment $xy$, and the velocity of the point $z$ with the segment $zy$, see Figure \ref{trailer}. The lengths of the two segments are fixed, and the length a horizontal curve is defined as the length of its $x$-projection to the plane. 
\begin{figure}[ht] 
\centering
\includegraphics[height=1.3 in]{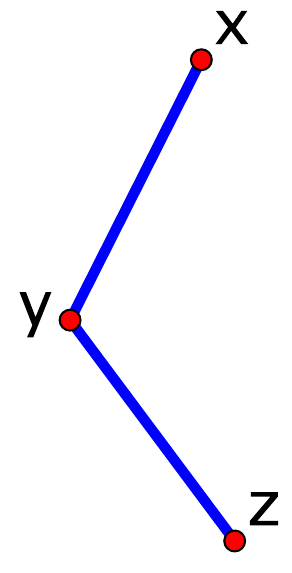}
\caption{Another 2-linkage.}	
\label{trailer}
\end{figure}

Rather than modeling bicycle, this is a model of tractor trailer, which is a well studied topic; see, e.g., \cite{TMS} and the references therein. 
Unfortunately, we did not obtain meaningful results in this case, but we think that it merits a study. See Figure \ref{prb} for an example of the planar projections of a geodesic.
\begin{figure}[ht] 
\centering
\includegraphics[height=1.75 in]{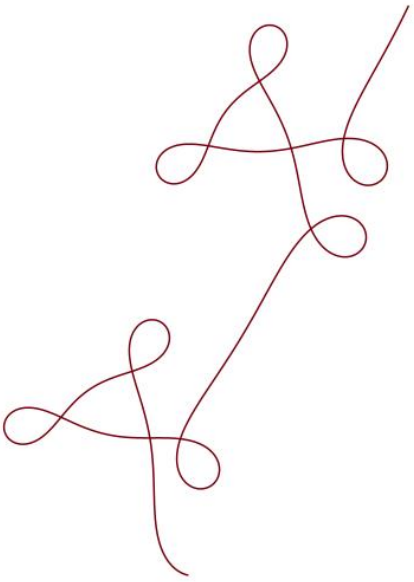}
\includegraphics[height=1.75 in]{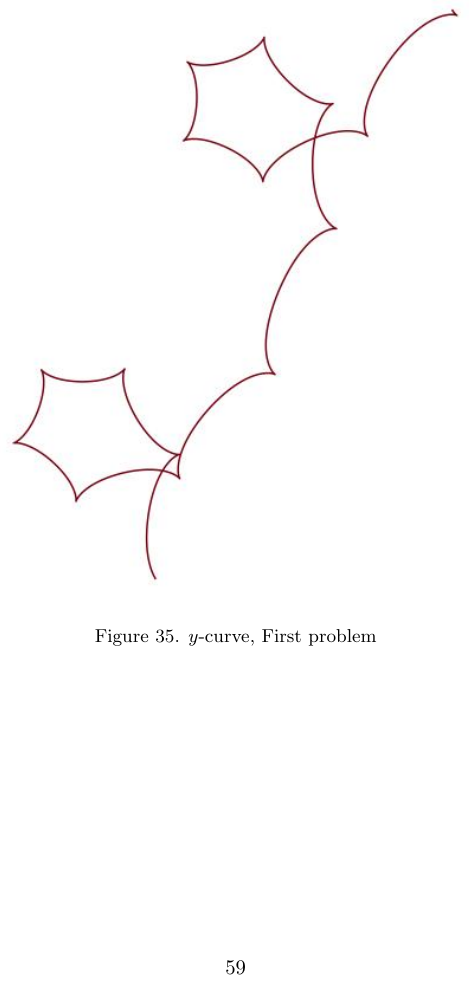}
\includegraphics[height=1.75 in]{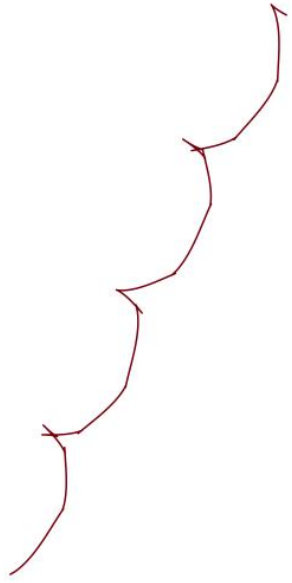}
\caption{The $x,y$, and $z$ projections of a geodesic.}	
\label{prb}
\end{figure}
\item Why stop at 2-linkages?  There are obvious $n$-linkage problems to consider, where the non-holonomic constraints apply to the midpoints of each segment, and one wants to extremize the length of the trajectory of either vertex (they have the same length).

In the case of three or more segments, the {\it topology} of the linkage is not unique.  In the case of 3-linkages, one can have a 
``linear" linkage, or a ``spokewheel" linkage, see Figure \ref{3link}.
\begin{figure}[ht] 
\centering
\includegraphics[width=.9\textwidth]{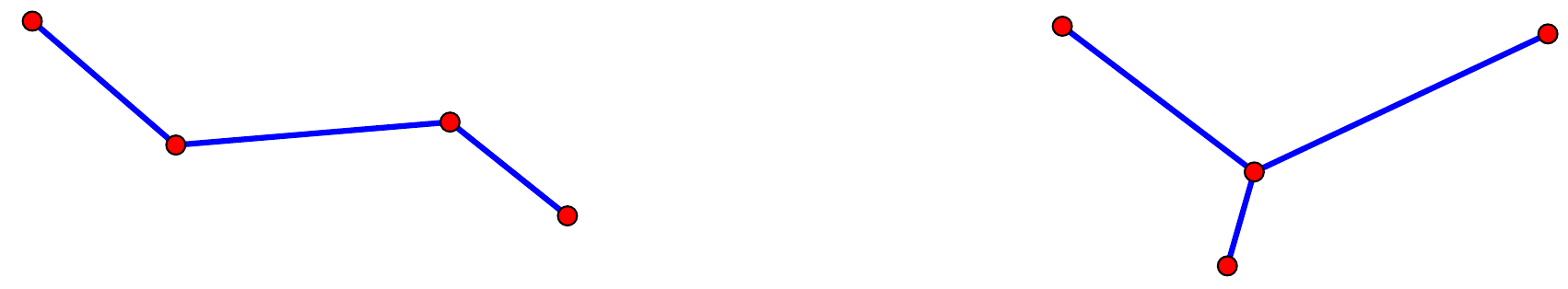}
\caption{Two topologically different 3-linkages.}	
\label{3link}
\end{figure}

Our last figures give examples of the trajectories in the ``linear" and ``spokewheel" cases.  Perhaps one or both
of these problems is integrable, as well extensions to linkages with four or more links. See Figure \ref{lin3diff},
\ref{wheel}. 
\end{itemize}

\end{document}